\documentclass[12pt]{article}
\usepackage{graphicx}
\usepackage{amsmath,amsthm,amssymb,enumerate}
\usepackage{euscript,mathrsfs}
\usepackage{color}
\usepackage{dsfont}
\usepackage[left=2cm,right=2cm,top=3.5cm,bottom=3.5cm]{geometry}
\usepackage{color}
\usepackage[framemethod=tikz]{mdframed}
\allowdisplaybreaks

\usepackage{soul}

\catcode`\@=11 \@addtoreset{equation}{section}

\catcode`\@=12

\newtheorem{Theorem}{Theorem}[section]
\newtheorem{Proposition}[Theorem]{Proposition}
\newtheorem{Lemma}[Theorem]{Lemma}
\newtheorem{Corollary}[Theorem]{Corollary}

\theoremstyle{definition}
\newtheorem{Definition}[Theorem]{Definition}

\newtheorem{Remark}[Theorem]{Remark}

\newcommand{\bTheorem}[1]{
	\begin{Theorem} \label{T#1} }
	\newcommand{\eT}{\end{Theorem}}

\newcommand{\bProposition}[1]{
	\begin{Proposition} \label{P#1}}
	\newcommand{\eP}{\end{Proposition}}

\newcommand{\bLemma}[1]{
	\begin{Lemma} \label{L#1} }
	\newcommand{\eL}{\end{Lemma}}

\newcommand{\bCorollary}[1]{
	\begin{Corollary} \label{C#1} }
	\newcommand{\eC}{\end{Corollary}}

\newcommand{\bRemark}[1]{
	\begin{Remark} \label{R#1} }
	\newcommand{\eR}{\end{Remark}}

\newcommand{\bDefinition}[1]{
	\begin{Definition} \label{D#1} }
	\newcommand{\eD}{\end{Definition}}

\newcommand{\tvB}{\widetilde{\vB}}

\newcommand{\Del}{\Delta_x}

\newcommand{\bB}{\vc{B}_B}

\newcommand{\vB}{\vc{B}}

\newcommand{\vuB}{\vc{u}_B}

\newcommand{\bfphi}{\boldsymbol{\varphi}}

\newcommand{\bFormula}[1]{
	\begin{equation} \label{#1}}
	\newcommand{\eF}{\end{equation}}

\newcommand{\Ov}[1]{\overline{#1}}

\newcommand{\Curl}{{\bf curl}_x}

\newcommand{\aleq}{\stackrel{<}{\sim}}

\newcommand{\ageq}{\stackrel{>}{\sim}}

\newcommand{\vr}{\varrho}

\newcommand{\tvt}{\tilde \vt}

\newcommand{\vt}{\vartheta}
\newcommand{\vu}{\vc{u}}
\newcommand{\vm}{\vc{m}}

\newcommand{\vc}[1]{{\bf #1}}

\newcommand{\Div}{{\rm div}_x}
\newcommand{\Grad}{\nabla_x}

\newcommand{\dx}{\,{\rm d} {x}}

\newcommand{\dt}{\,{\rm d} t }

\newcommand{\vU}{\vc{U}}

\newcommand{\intO}[1]{\int_{\Omega} #1 \ \dx}

\newcommand{\vv}{\vc{v}}

\newcommand{\D}{{\rm d}}

\newcommand{\ep}{\varepsilon}

\newcommand{\vtB}{\vt_B}

\newcommand{\br}{ \nonumber \\ }

\def\softd{{\leavevmode\setbox1=\hbox{d}%
		\hbox to 1.05\wd1{d\kern-0.4ex{\char039}\hss}}}
\definecolor{Cgrey}{rgb}{0.85,0.85,0.85}
\definecolor{Cblue}{rgb}{0.50,0.85,0.85}
\definecolor{Cred}{rgb}{1,0,0}
\definecolor{fancy}{rgb}{0.10,0.85,0.10}

\newcommand\Cbox[2]{%
	\newbox\contentbox%
	\newbox\bkgdbox%
	\setbox\contentbox\hbox to \hsize{%
		\vtop{
			\kern\columnsep
			\hbox to \hsize{%
				\kern\columnsep%
				\advance\hsize by -2\columnsep%
				\setlength{\textwidth}{\hsize}%
				\vbox{
					\parskip=\baselineskip
					\parindent=0bp
					#2
				}%
				\kern\columnsep%
			}%
			\kern\columnsep%
		}%
	}%
	\setbox\bkgdbox\vbox{
		\color{#1}
		\hrule width  \wd\contentbox %
		height \ht\contentbox %
		depth  \dp\contentbox
		\color{black}
	}%
	\wd\bkgdbox=0bp%
	\vbox{\hbox to \hsize{\box\bkgdbox\box\contentbox}}%
	\vskip\baselineskip%
}

\mdfdefinestyle{MyFrame}{%
	linecolor=black,
	outerlinewidth=1pt,
	roundcorner=5pt,
	innertopmargin=\baselineskip,
	innerbottommargin=\baselineskip,
	innerrightmargin=10pt,
	innerleftmargin=10pt,
	backgroundcolor=white!20!white}

\begin{document}


\title{Mathematical theory of compressible magnetohydrodynamics  
driven by non--conservative boundary conditions}

\author{Eduard Feireisl
		\thanks{The work of E.F. was partially supported by the
			Czech Sciences Foundation (GA\v CR), Grant Agreement
			21--02411S. The Institute of Mathematics of the Academy of Sciences of
			the Czech Republic is supported by RVO:67985840. This work was partially supported by the Thematic Research Programme, University of Warsaw, Excellence Initiative Research University.} \and Piotr Gwiazda
		\and \and Young--Sam Kwon
		\thanks{The work of Y.--S. Kwon was partially supported by
			the National Research Foundation of Korea (NRF2022R1F1A1073801)}
		\and Agnieszka \'Swierczewska-Gwiazda
		\thanks{The work of P. G. and A. \'S-G. was partially supported by  National Science Centre
			(Poland),  agreement no   2021/43/B/ST1/02851.}}

\date{}

\maketitle

\medskip

\centerline{Institute of Mathematics of the Academy of Sciences of the Czech Republic}

\centerline{\v Zitn\' a 25, CZ-115 67 Praha 1, Czech Republic}

\medskip

\centerline{Institute of Mathematics of Polish Academy of Sciences}
\centerline{\'Sniadeckich 8, 00-956 Warszawa, Poland}

\medskip

\centerline{Department of Mathematics, Dong-A University}

\centerline{Busan 49315, Republic of Korea}

\medskip

\centerline{Institute of Applied Mathematics and Mechanics, University of Warsaw}

\centerline{Banacha 2, 02-097 Warsaw, Poland}

\begin{abstract}
	
We propose a new concept of weak solution to the equations of compressible magnetohydrodynamics driven
by large boundary data. The system of the underlying field equations is solvable globally in time
in the out of equilibrium regime characteristic for turbulence. The weak solutions comply with the
weak--strong uniqueness principle; they coincide with the classical solution of the problem as long as the latter exists. The choice of constitutive relations is motivated by applications in stellar magnetoconvection.

\end{abstract}



\tableofcontents

\section{Introduction}
\label{i}

The large scale dynamics arising in stellar magnetoconvection is driven in an essential way by the boundary
conditions, see Gough \cite{Gough}. The flux separation and pattern formation arise in the turbulent regime when the fluid flow is far from equilibrium. As the underlying field equations are non--linear, the existence of global in time smooth
solutions in such a regime is not known. Indeed the vast majority of \emph{global} existence results in the class of classical solutions
is restricted to the initial state close to a stable equilibrium, see a.g. Matsumura and Nishida \cite{MANI},
Valli \cite{VAL}, among others. What is more, the recent results of Merle et al.
\cite{MeRaRoSz}, \cite{MeRaRoSzbis} indicate that classical solutions may develop singularities in a finite time.

In view of the above arguments, the concept of \emph{weak solution} represents a suitable alternative to restore
global--in--time existence even for problems with large data and solutions remaining out of equilibrium in the long run. We consider a mathematical model of fully compressible
three dimensional fluid convection driven by an externally imposed magnetic field. The fluid occupies a bounded
domain $\Omega \subset R^3$ with regular boundary. The time evolution of the density $\vr = \vr(t,x)$, the
(absolute) temperature $\vt = \vt(t,x)$, the velocity field $\vu = \vu(t,x)$, and the magnetic field
$\vB = \vB(t,x)$ is governed by the \emph{compressible MHD system} of field equations (cf.~\cite{DUFE2}):

\begin{mdframed}[style=MyFrame]
	
	{\bf Equation of continuity}
	
	\begin{equation} \label{p1}
		\partial_t \vr + \Div (\vr \vu) = 0.
	\end{equation}
	
	\noindent	{\bf Momentum equation}
	
	\begin{equation} \label{p2}
		\partial_t (\vr \vu) + \Div (\vr \vu \otimes \vu) + \Grad p (\vr, \vt) = \Div \mathbb{S}(\vt, \Grad \vu) + \Curl  \vc{B} \times \vc{B} +   \vr \vc{g}.
	\end{equation}

\noindent
{\bf Induction equation}

\begin{equation} \label{p3}
	\partial_t \vc{B} + \Curl (\vc{B} \times \vu ) + \Curl (\zeta (\vt) \Curl  \vc{B} ) = 0,\
	\Div \vc{B} = 0.
	\end{equation}

	\noindent {\bf Internal energy balance}
	\begin{align}
		\partial_t ( \vr e(\vr, \vt) ) + \Div (\vr e(\vr, \vt) \vu) &+ \Div \vc{q}(\vt, \Grad \vt) \br &=
		\mathbb{S}(\vt, \Grad \vu):\Grad \vu + \zeta(\vt) |\Curl \vB|^2  - p(\vr, \vt) \Div \vu.
		\label{p4}
	\end{align}	

\end{mdframed}

The pressure $p = p(\vr, \vt)$ and the internal energy $e=e(\vr, \vt)$ are interrelated through \emph{Gibbs' equation}
\begin{equation} \label{w1}
	\vt D s = De + p D \left( \frac{1}{\vr} \right),
\end{equation}
where $s = s(\vr, \vt)$ is the entropy. Consequently, the internal energy balance \eqref{p4} may be reformulated
in the form of \emph{entropy balance equation}
\begin{align}
	\partial_t (\vr s(\vr, \vt)) + \Div (\vr s(\vr, \vt) \vu) &+ \Div \left( \frac{\vc{q}(\vt, \Grad \vt)}{\vt} \right) \br &=
	\frac{1}{\vt} \left( \mathbb{S}(\vt, \Grad \vu) : \Grad \vu - \frac{\vc{q}(\vt, \Grad \vt) \cdot \Grad \vt}{\vt} + \zeta(\vt) |\Curl \vB |^2 \right), \label{w2}
\end{align}
see e.g. Gough \cite{Gough}, Weiss and Proctor \cite{Weiss}, Tao et al. \cite{Taoetal}.

The right--hand side of equation \eqref{p4} represents the entropy production rate, which, in accordance with the Second law of thermodynamics, must be non--negative. Accordingly, we consider \emph{Newtonian fluid}, 
with the viscous stress tensor
\begin{equation} \label{c7}
	\mathbb{S}(\vt, \Grad \vu) = \mu (\vt) \left( \Grad \vu + \Grad^t \vu - \frac{2}{3} \Div \vu \mathbb{I} \right) +
	\eta(\vt) \Div \vu \mathbb{I},
\end{equation}
where the viscosity coefficients $\mu > 0$ and $\eta \geq 0$ are continuously differentiable functions of the temperature. Similarly, the heat flux obeys \emph{Fourier's law},
\begin{equation} \label{c9}
	\vc{q}(\vt, \Grad \vt)= - \kappa (\vt) \Grad \vt,
\end{equation}	
where the heat conductivity coefficient $\kappa > 0$ is a continuously differentiable function of the temperature.

The existence of global--in--time weak solutions for the compressible MHD system \eqref{p1}--\eqref{c9}
was shown in \cite{DUFE2} under certain physically relevant restrictions imposed on constitutive equations and
transport coefficients. The boundary conditions considered in \cite{DUFE2} are conservative,
\begin{equation} \label{conbc}
\vu|_{\partial \Omega} = 0, \ \Grad \vt \cdot \vc{n}|_{\partial \Omega} = 0,
\ \vB \cdot \vc{n}|_{\partial \Omega} = 0,\ \left[ \vB \times \vu + \zeta \Curl \vB \right] \times \vc{n}|_{\partial \Omega} = 0
\end{equation}
characteristic for closed systems. If the driving force $\vc{g}$ in the momentum equation is conservative,
meaning $\vc{g} = \Grad G$, $G = G(x)$, the total energy of the system is conserved and the dynamics obeys the rather ``boring'' scenario formulated by the celebrated statement of Clausius:
\begin{quotation}
	
	\emph{
		The energy of the world is constant; its entropy tends to a maximum.}
	
\end{quotation}
A rigorous mathematical proof of this statement for the Navier--Stokes--Fourier system was given in
\cite{FP20}.

As pointed out by Gough \cite{Gough}, a rich fluid behaviour is conditioned by a proper choice of boundary
conditions. Motivated by models in astrophysics, we suppose the boundary $\partial \Omega$ is impermeable
and the tangential component of the normal viscous stress vanishes on it,
	\begin{equation} \label{p5}
	\vu \cdot \vc{n}|_{\partial \Omega} = 0,\ ( \mathbb{S} \cdot \vc{n} ) \times \vc{n}|_{\partial \Omega} = 0.
\end{equation}
However, the theory presented in this paper can accommodate more general boundary conditions for the velocity 
field discussed in Section \ref{C}.  

Similarly to the well known Rayleigh--B\' enard problem, see e.g. Davidson \cite{DAVI}, we impose the inhomogeneous
Dirichlet boundary conditions for the temperature,
\begin{equation} \label{p7}
	\vt|_{\partial \Omega} = \vtB ,\ \vtB > 0.
\end{equation}

To incorporate the effect of an exterior magnetic field, we suppose that $\vc{B}$ is (not necessarily small)
perturbation of a background magnetic field $\bB$, $\Div \bB = 0$,
\[
\vB = \bB + \vc{b}.
\]
In the context of astrophysics, we may think of $\bB$ as being the magnetic field imposed by a massive star. As for $\vc{b}$,  we
suppose either
\begin{equation} \label{p8}
	\vc{b} \times \vc{n}|_{\partial \Omega} = 0 \ \Rightarrow \ \vc{B} \times \vc{n}|_{\partial \Omega} = \bB
	\times \vc{n}|_{\partial \Omega},
	\end{equation}
or
\begin{equation} \label{p9}
	\vc{b} \cdot \vc{n}|_{\partial \Omega} = 0 \ \Rightarrow \ \vc{B} \cdot \vc{n}|_{\partial \Omega} = 	\bB
	\cdot \vc{n}|_{\partial \Omega}.
	\end{equation}
The boundary condition \eqref{p8} is of Dirichlet type, compatible with the requirement of solenoidality of $\vB$. 
The condition \eqref{p9} is of flux type, similar to the first condition \eqref{p5} and must be accompanied 
by another flux condition related to the electric field, namely
\begin{equation} \label{w4a}
	\Big[  \vB \times \vu + \zeta  \Curl \vB \Big] \times \vc{n}|_{\partial \Omega} = 0.	
\end{equation}
The flux in \eqref{w4a} may be inhomogeneous as well. The borderline case $\vc{B}_B = 0$ in  \eqref{p8} and 
\eqref{p9} corresponds to a perfectly isolating and perfectly conducting boundary, respectively (see e.g. 
Alekseev \cite{Aleks}). 

For possibly technical but so far unsurmountable reasons, a mathematically tractable \emph{weak formulation} of compresible MHD system cannot be based on
merely rewriting the system \eqref{p1}--\eqref{p4} in the sense of distributions. The available {\it a priori} bounds are not strong enough to render certain terms, notably $p(\vr, \vt) \Div \vu$ in
the internal energy balance \eqref{p4}, integrable. A remedy proposed in \cite{DUFE2} is replacing
\eqref{p4} by the entropy equation \eqref{w2}. Unfortunately, the entropy production rate is {\it a priori}
bounded only in the non--reflexive space $L^1$ of integrable functions. As a result, any approximation scheme
provides merely an \emph{inequality}
\begin{align}
	\partial_t (\vr s(\vr, \vt)) + \Div (\vr s(\vr, \vt) \vu) &+ \Div \left( \frac{\vc{q}(\vt, \Grad \vt)}{\vt} \right) \br &\geq
	\frac{1}{\vt} \left( \mathbb{S}(\vt, \Grad \vu) : \Grad \vu - \frac{\vc{q}(\vt, \Grad \vt) \cdot \Grad \vt}{\vt} + \zeta(\vt) |\Curl \vB |^2 \right) \label{w2a}
\end{align}
satisfied in the sense of distributions.

Replacing entropy equation by inequality definitely enlarges the class of possible
solutions and gives rise to an underdetermined problem. To save well posedness, at least formally, the integrated \emph{total}
energy balance/inequality is appended to the system \eqref{p1}--\eqref{p3}, \eqref{w2a} in \cite{DUFE2}.
The total energy balance reads
\begin{align}
\partial_t &\left( \frac{1}{2} \vr |\vu|^2 + \vr e(\vr, \vt) + \frac{1}{2} |\vB |^2 \right) +
\Div \left[ \left( \frac{1}{2} \vr |\vu|^2 + \vr e(\vr, \vt) + \frac{1}{2} |\vB |^2 + p(\vr, \vt) \right) \vu \right] \br
&+ \Div \Big[ (\vB \times \vu) \times \vB + \zeta(\vt) \Curl \vB \times \vB \Big]  -
\Div( \mathbb{S} (\vt, \Grad \vu)\cdot \vu ) + \Div \vc{q}(\vt, \Grad \vt) = \vr \vc{g} \cdot \vu.	
	\label{teb}
	\end{align}
If the system is \emph{conservative}, the total energy flux vanishes on the boundary and \eqref{teb} integrated over $\Omega$ yields
\begin{equation} \label{teb1}
\frac{\D}{\dt} \intO{ \left( \frac{1}{2} \vr |\vu|^2 + \vr e(\vr, \vt) + \frac{1}{2} |\vB |^2 \right)  } =
\intO{ \vr \vc{g} \cdot \vu }.	
	\end{equation}
It turns out that the equations \eqref{p1}--\eqref{p3}, the inequality \eqref{w2a}, together with the integral identity \eqref{teb1}
represent a suitable weak formulation for the \emph{conservative} system. Indeed the existence of global--in--time weak solutions for any
finite energy (initial) data was proved in \cite{DUFE2}. The weak solutions comply with a natural \emph{compatibility} principle. Any sufficiently smooth weak solution is a classical solution of the problem. In addition, the weak solutions satisfy the weak--strong uniqueness principle: A weak solution coincides with
the strong solution on the life span of the latter, see \cite{FeiNov10}.

The situation becomes more delicate for open (non--conservative) systems, with inhomogeneous Dirichlet boundary conditions. The heat flux $\vc{q} \cdot \vc{n}$ as well as the induction flux (the normal component of the electric field) $[ \vB \times \vu \times \vB + \zeta \Curl \vB \times \vB ] \cdot \vc{n}$ 
do not necessarily 
vanish on $\partial \Omega$ and give rise to additional uncontrollable sources of energy in \eqref{teb1}. In \cite{ChauFei}
(see also the monograph \cite{FeiNovOpen}), a new approach has been developed to handle the inhomogeneous Dirichlet
boundary conditions for the temperature. The total energy balance \eqref{teb1} is replaced by a similar inequality for
the ballistic energy
\[
\left[ \frac{1}{2} \vr |\vu|^2 + \vr e(\vr, \vt) + \frac{1}{2} |\vB |^2 - \tvt \vr s (\vr, \vt) \right],
\]
where $\tvt$ is an arbitrary extension of the boundary temperature $\vtB$ inside $\Omega$. Pursuing the same
strategy, we introduce a ``magnetic'' variant of the ballistic energy
\begin{equation} \label{w3}
E_{BM} = \left[	\frac{1}{2} \vr |\vu|^2 + \vr e + \frac{1}{2} |\vB |^2 - \tvt \vr s - \bB \cdot \vB \right],
\end{equation}
After a straightforward manipulation, we deduce equality
\begin{align}
	\frac{\D }{\dt} &\intO{ \left( \frac{1}{2} \vr |\vu|^2 + \vr e + \frac{1}{2} |\vB |^2 - \tvt \vr s - \bB \cdot \vB \right) } \br
	&+ \intO{ 	\frac{\tvt}{\vt} \left( \mathbb{S} : \Grad \vu - \frac{\vc{q} \cdot \Grad \vt}{\vt} + \zeta | \Curl \vB |^2 \right) } \br
	&= - \intO{\left( \vr s \partial_t \tvt + \vr s \vu \cdot \Grad \tvt + \frac{\vc{q}}{ \vt} \cdot \Grad \tvt                    \right)      } \br
	&\quad -  \intO{  \Big( \vB \cdot \partial_t \bB - (\vB \times \vu) \cdot \Curl \bB -
		\zeta \Curl \ \vB \cdot \Curl \bB \Big)    } \br
	&\quad + \intO{ \vr \vc{g} \cdot \vu }
	\label{w4}
\end{align}
as soon as $\tvt|_{\partial \Omega} = \vtB$, $\tvt > 0$, $\vB \times \vc{n}|_{\partial \Omega} =
\bB \times \vc{n}|_{\partial \Omega}$. The same is true if $\vB \cdot \vc{n}|_{\partial \Omega} =
\bB \cdot \vc{n}|_{\partial \Omega}$ and the vanishing tangential flux condition \eqref{w4a} holds.  
Note carefully that \eqref{w4} does not contain the boundary flux integrals; whence it is suitable for a weak formulation of the problem. 

The mathematical theory we propose is based on imposing the weak form of the equations \eqref{p1}--\eqref{p3}, together with
the entropy inequality \eqref{w2a}, and the ballistic energy balance \eqref{w4} as a weak formulation of
the compressible MHD system. The basic hypotheses concerning the state equation and the transport coefficients
as well as the exact definition of weak solution are given in Section \ref{w}. In Section \ref{r} we introduce the
basic tool to investigate stability properties in the class of weak solutions -- the \emph{relative energy inequality}. In Section \ref{ws}, we show that the weak solutions enjoy the weak--strong uniqueness property.
In Section \ref{e}, the existence of global--in--time weak solutions is established by means of a multilevel approximation
scheme. The paper is concluded in Section \ref{C} by a brief discussion on possible extensions of the theory to a larger class of boundary conditions.

\section{Weak formulation}
\label{w}

We start by a list of structural hypotheses imposed on the equation of state and the transport coefficients.

\subsection{Equation of state}

Our choice of the equation of state is motivated by \cite[Chapter 4]{FeiNovOpen}. In particular, we consider the
radiation contribution to the pressure/internal energy relevant to problems in astrophysics, cf.
Battaner \cite{BATT}. In accordance with Gibbs' relation \eqref{w1} we suppose that
\begin{align}
	p(\vr, \vt) &= p_M (\vr, \vt) + p_R (\vt), \ p_M(\vr,\vt) = \vt^{\frac{5}{2}} P \left( \frac{\vr}{\vt^{\frac{3}{2}}  } \right), \ p_R (\vt) = \frac{a}{3} \vt^4,\br
	e(\vr, \vt) &= e_M (\vr, \vt) + e_R (\vr, \vt),\ e_M(\vr,\vt) =  \frac{3}{2} \frac{\vt^{\frac{5}{2}} }{\vr} P \left( \frac{\vr}{\vt^{\frac{3}{2}}  } \right),\ e_R (\vr ,\vt) = \frac{a}{\vr} \vt^4, \ a > 0,
	\label{c1}	
\end{align}
where the function $P \in C^1[0,\infty)$ satisfies
\begin{equation} \label{c2}
	P(0) = 0,\ P'(Z) > 0 \ \mbox{for}\ Z \geq 0,\ 0 < \frac{ \frac{5}{3} P(Z) - P'(Z) Z }{Z} \leq c \ \mbox{for}\ Z > 0.
\end{equation} 	
This implies, in particular, that $Z \mapsto P(Z)/ Z^{\frac{5}{3}}$ is decreasing, and we suppose
\begin{equation} \label{c3}
	\lim_{Z \to \infty} \frac{ P(Z) }{Z^{\frac{5}{3}}} = p_\infty > 0.
\end{equation}
Note that $p_M$, $e_M$ satisfy the relation characteristic for monoatomic gas,
\[
p_M = \frac{2}{3} \vr e_M,
\]
$p_R$ is the radiation pressure with the associated internal energy $e_R$, and \eqref{c3} corresponds to the
presence of electron pressure under the degenerate gas regime.
It follows from \eqref{c2} that $p$ and $e$ satisfy the \emph{hypothesis of thermodynamic stability}:
\begin{equation} \label{c4}
	\frac{ \partial p (\vr, \vt) }{\partial \vr} > 0,\
	\frac{ \partial e (\vr, \vt) }{\partial \vt} > 0.
\end{equation}

In accordance with \eqref{w1}, the entropy takes the form
\begin{equation} \label{c5}
	s(\vr, \vt) = s_M(\vr, \vt) + s_R (\vr, \vt),\ s_M(\vr, \vt) = \mathcal{S} \left( \frac{\vr}{\vt^{\frac{3}{2}} } \right),\ s_R(\vr, \vt) = \frac{4a}{3} \frac{\vt^3}{\vr},
\end{equation}
where
\begin{equation} \label{c6}
	\mathcal{S}'(Z) = -\frac{3}{2} \frac{ \frac{5}{3} P(Z) - P'(Z) Z }{Z^2}.
\end{equation}	

\subsection{Transport coefficients}

We suppose the viscosity coefficients in the Newtonian stress $\mathbb{S}(\vt, \Grad \vu)$ are continuously
differentiable functions of the temperature satisfying
\begin{align}
	0 < \underline{\mu} \left(1 + \vt^\alpha \right) &\leq \mu(\vt) \leq \Ov{\mu} \left( 1 + \vt^\alpha \right),\
	|\mu'(\vt)| \leq c \ \mbox{for all}\ \vt \geq 0, \br
	0 &\leq  \eta(\vt) \leq \Ov{\eta} \left( 1 + \vt^\alpha \right) \ \mbox{for some}\  \frac{1}{2} \leq \alpha \leq 1.
	\label{c8}
\end{align}
Note that $\alpha = \frac{1}{2}$ corresponds to Sutherland's law relevant in astrophysics, 
see Yang et al. \cite{YangWMC}.

Similarly, the heat conductivity coefficient in the Fourier heat flux $\vc{q}(\vt, \Grad \vt)$ is a
continuously differentiable function of the temperature satisfying
\begin{equation} \label{c10}
	0 < \underline{\kappa} \left(1 + \vt^\beta \right) \leq  \kappa(\vt) \leq \Ov{\kappa} \left( 1 + \vt^\beta \right) \ \mbox{for some}\  \beta \geq 3.
\end{equation}
Here, the case $\beta = 3$ reflects the effect of radiation.

Finally, we suppose the magnetic diffusivity coefficient $\zeta = \zeta(\vt)$ is a continuously differentiable function of the temperature,
\begin{equation} \label{c11}
	0 < \underline{\zeta}(1 + \vt) \leq \zeta (\vt) \leq \Ov{\zeta}(1 + \vt),\ |\zeta' (\vt) | \leq c
	\ \mbox{for all}\ \vt \geq 0.	
\end{equation}

\subsection{Boundary data}

We suppose the background magnetic field $\bB = \bB(t,x)$ is sufficiently
smooth at least of class $C^{1,2}([0,T] \times \Ov{\Omega}; R^3)$,
\begin{equation} \label{BC1}
	\Div \bB = 0 \ \mbox{in}\ (0,T) \times \Omega.
	\end{equation}
Similarly, we suppose the boundary temperature $\vtB$ can be extended inside $\Omega$,
$\vtB \in C^{1,2}([0,T] \times \Ov{\Omega} )$,
\begin{equation} \label{BC2}
	\vtB > 0 \ \mbox{in}\ [0,T] \times \Ov{\Omega}.
	\end{equation}

\subsection{Weak solutions}

Before discussing the concept of weak solution, it is useful to introduce the function spaces
\begin{align}
H_{0, \tau}(\Omega;R^3) &= \left\{ \vc{b} \in L^2(\Omega; R^3) \ \Big| \
\Curl \vc{b} \in L^2(\Omega; R^d),\ \Div \vc{b} = 0,\ \vc{b} \times \vc{n}|_{\partial \Omega} = 0 \right\},\br
H_{0, n}(\Omega;R^3) &= \left\{ \vc{b} \in L^2(\Omega; R^3) \ \Big| \
\Curl \vc{b} \in L^2(\Omega; R^d),\ \Div \vc{b} = 0,\ \vc{b} \cdot \vc{n}|_{\partial \Omega} = 0 \right\}.
	\nonumber
	\end{align}
Both $H_{0,\tau}$ and $H_{0,n}$ are endowed with the Hilbert norm \[ 
\| \vc{b} \|^2_{H_0} = \| \Curl \vc{b} \|^2_{L^2(\Omega; R^3)} + \| \vc{b} \|^2_{L^2(\Omega; R^3)},
\]
see e.g. Alexander and Auchmuty 
\cite{AuAl}.

\begin{mdframed}[style=MyFrame]
	
	\begin{Definition} \label{De1} {\bf (Weak solution to the compressible MHD system)}

A quantity $(\vr, \vt, \vu, \vB)$ is termed \emph{weak solution} of the compressible MHD system
\eqref{p1}--\eqref{p4}, with the boundary conditions \eqref{p5}, \eqref{p7}, \eqref{p8} (respectively
\eqref{p9}), and the initial data
\begin{equation} \label{ID}
	\vr(0, \cdot) = \vr_0,\ \vt(0, \cdot) = \vt_0,\ \vu(0, \cdot) = \vu_0,\ \vB(0,\cdot) = \vB_0,
	\end{equation}
if the following holds:
\begin{itemize}
	\item{\bf Equation of continuity.}
	$\vr \in C_{\rm weak}([0,T]; L^{\frac{5}{3}} (\Omega))$, $\vr \geq 0$ in $(0,T) \times \Omega$, and the integral
	identity
\begin{equation} \label{w5}
	\int_0^\tau \intO{ \Big( \vr \partial_t \varphi + \vr \vu \cdot \Grad \varphi \Big) } \ dt = -
	\left[ \intO{ \vr \varphi } \right]_{t=0}^{t = \tau},\ \vr(0, \cdot) = \vr_0,
\end{equation}
holds for any $\varphi \in C^1([0,T] \times \Ov{\Omega})$. 	
In addition, the renormalized version of \eqref{w5}
\begin{equation} \label{w7}
	\int_0^\tau \intO{ \left( b(\vr) \partial_t \varphi + b(\vr) \vu \cdot \Grad \varphi +
		\Big( b(\vr) - b'(\vr) \vr \Big) \Div \vu \right) } \dt = \left[ \intO{ b(\vr) \varphi  } \right]_{t = 0}^{t= \tau}
\end{equation}
holds for any $\varphi \in C^1([0,T] \times \Ov{\Omega})$ and any $b \in C^1(R)$, $b' \in C_c(R)$.

\item{\bf Momentum equation.}

$\vr \vu \in C_{\rm weak}([0,T]; L^{\frac{5}{4}}(\Omega; R^3)$,
$\vu \in L^2(0,T; W^{1,q}(\Omega; R^3))$, $q = \frac{8}{5 - \alpha}$, $\vu \cdot \vc{n}|_{\partial \Omega} = 0$,
and the integral identity
\begin{align}
	\int_0^\tau &\intO{ \Big( \vr \vu \cdot \partial_t \bfphi + \vr \vu \otimes \vu : \Grad \bfphi +
		p(\vr, \vt) \Div \bfphi \Big) } \br &= \int_0^\tau \intO{ \mathbb{S}(\vt, \Grad \vu) : \Grad \bfphi } \dt -
	\int_0^\tau \intO{ \left( \vB \otimes \vB - \frac{1}{2} |\vB |^2 \mathbb{I} \right) : \Grad \bfphi } \dt \br &-
	\int_0^\tau \intO{ \vr \vc{g} \cdot \bfphi } \dt  + \left[  \intO{ \vr \vu \cdot \bfphi }
	\right]_{t = 0}^{t = \tau},\ (\vr \vu)(0, \cdot) = \vr_0 \vu_0	
	\label{w9}
\end{align}
holds for any $\bfphi \in C^1([0,T] \times \Ov{\Omega}; R^3)$, $\bfphi \cdot \vc{n}|_{\partial \Omega} = 0$.

\item {\bf Induction equation.}
$\vB \in C_{\rm weak}([0,T]; L^2(\Omega; R^3))$,
\begin{equation} \label{w10}
	\Div \vc{B}(\tau, \cdot) = 0 \ \mbox{for any}\ \tau \in [0,T].
	\end{equation}
In addition
\begin{equation} \label{w10a}
	(\vB - \bB) \in L^2(0, T; H_{0, \tau}(\Omega; R^3))
	\end{equation}
if the boundary condition \eqref{p8} is imposed or
\begin{equation} \label{w10b}
	(\vB - \bB) \in L^2(0, T; H_{0, n}(\Omega; R^3))
\end{equation}
for the boundary conditions \eqref{p9}, \eqref{w4a}.
The integral identity
\begin{align}
	\int_0^\tau &\intO{ \Big( \vB \cdot \partial_t \bfphi - (\vB \times \vu) \cdot \Curl \bfphi -
		\zeta (\vt) \Curl \ \vB \cdot \Curl \ \bfphi \Big) } \dt \br &= \left[ \intO{ \vB \cdot \bfphi  } \right]_{t=0}^{t = \tau},\ \vB(0, \cdot) = \vB_0,
	\label{w11}
\end{align}	
holds for any $\bfphi \in C^1([0,T] \times \Ov{\Omega}; R^3)$, $\Div \bfphi = 0$,
\begin{equation} \label{w11a}
\bfphi \times \vc{n}|_{\partial \Omega} = 0
\end{equation}
in the case of boundary conditions \eqref{p8},
\begin{equation} \label{w11b}
	\bfphi \cdot \vc{n}|_{\partial \Omega} = 0
\end{equation}
for the boundary conditions \eqref{p9}, \eqref{w4a}.

\item {\bf Entropy inequality.}
$\vt \in L^\infty(0,T; L^4(\Omega))$, $\vt > 0$ a.a. in $(0,T) \times \Omega$,
$\vt - \vtB \in L^2(0,T; W^{1,2}_0 (\Omega))$, $\log (\vt) \in L^2(0,T; W^{1,2}(\Omega))$, and the integral inequality
\begin{align}
	\int_{\tau_1}^{\tau_2} &\intO{ \left( \vr s (\vr, \vt) \partial_t \varphi + \vr s (\vr, \vt) \vu \cdot \Grad \varphi + \frac{\vc{q} (\vt, \Grad \vt) }{\vt} \cdot
		\Grad \varphi \right) } \dt \br & \leq  - \int_{\tau_1}^{\tau_2} \intO{ \frac{1}{\vt}
		\left( \mathbb{S}(\vt, \Grad \vu) : \Grad \vu - \frac{\vc{q} (\vt, \Grad \vt) \cdot \Grad \vt}{\vt} + \zeta (\vt) |\Curl \vB|^2 \right) \varphi }
	\dt\br & + \left[ \intO{ \vr s(\vr, \vt) \varphi } \right]_{t = \tau_1-}^{t = \tau_2+},\
	0 \leq \tau_1 \leq \tau_2 \leq T,\ \vr s(\vr, \vt)(0-, \cdot) = \vr_0 s(\vr_0, \vt_0), 	
	\label{w13}	
\end{align}
holds for any $\varphi \in C^1_c([0,T] \times \Omega)$, $\varphi \geq 0$.

\item {\bf Ballistic energy inequality.}
The inequality
\begin{align}
	&\intO{ \left( \frac{1}{2} \vr |\vu|^2 + \vr e(\vr, \vt) + \frac{1}{2} |\vB |^2 - \tvt \vr s (\vr, \vt) - \bB \cdot \vB \right) (\tau, \cdot) }
	\br&	+ \int_0^\tau \intO{ 	\frac{\tvt}{\vt} \left( \mathbb{S}(\vt, \Grad \vu) : \Grad \vu - \frac{\vc{q} (\vt, \Grad \vt) \cdot \Grad \vt}{\vt} + \zeta (\vt) | \Curl \vB |^2 \right) } \dt \br
	&\quad \leq - \int_0^\tau \intO{\left( \vr s (\vr, \vt) \partial_t \tvt + \vr s (\vr, \vt) \vu \cdot \Grad \tvt + \frac{\vc{q}(\vt, \Grad \vt)}{ \vt} \cdot \Grad \tvt                    \right)      } \br
	&\quad \quad	 -  \int_0^\tau \intO{  \Big( \vB \cdot \partial_t \bB - (\vB \times \vu) \cdot \Curl \bB -
		\zeta (\vt) \Curl \ \vB \cdot \Curl \bB \Big)    } \dt \br
	&\quad \quad	 + \int_0^\tau\intO{ \vr \vc{g} \cdot \vu } \dt \br
	&\quad \quad + \intO{ \left( \frac{1}{2} \vr_0 |\vu_0|^2 + \vr_0 e(\vr_0, \vt_0) + \frac{1}{2} |\vB_0 |^2 - \tvt (0, \cdot) \vr_0 s
		(\vr_0, \vt_0)  - \bB (0, \cdot) \cdot \vB_0 \right) }
	\label{w14}
	\end{align}
	holds for any  $\tvt \in C^1([0,T] \times \Ov{\Omega})$, $\tvt > 0$, $\tvt|_{\partial \Omega} = \vtB$.

	\end{itemize}		
		
		\end{Definition}

	\end{mdframed}

\bigskip

\begin{Remark} \label{Rw1}
	
	As the magnetic field $\vB$ satisfies \eqref{w10}, validity of
	\eqref{w11} can be extended to any (not necessarily solenoidal) test function $\bfphi$ satisfying
	\eqref{w11a}, \eqref{w11b}, respectively. Indeed, if $\bfphi$ satisfies the boundary condition \eqref{w11a},
	 we consider the potential $\Psi$ solving
	 \[
	 \Del \Psi = \Div \bfphi \ \mbox{in}\ \Omega,\ \Psi|_{\partial \Omega} = 0,
	 \]
	 and write
	 \[
	 \bfphi = (\bfphi - \Grad \Psi) + \Grad \Psi.
	 \]
	 On the one hand, as $\Psi$ vanishes on $\partial \Omega$ and $\vB$ is solenoidal, we have
	 \[
	 \intO{ \vc{B} \cdot \Grad \Psi } = \intO{ \vc{B} \cdot \partial_t \Grad \Psi} = 0;
	 \]
	 whence the integral identity \eqref{w11} holds for $\Grad \Psi$. On the other hand, the function
	 $(\bfphi - \Grad \Psi)$ is solenoidal, and, as $\Grad \Psi \times \vc{n}|_{\partial \Omega} = 0$,
	 a legal test function in \eqref{w11}.
	
	 If $\bfphi$ satisfies \eqref{w11b}, we simply consider
	 \[
	 \Del \Psi = \Div \bfphi,\ \Grad \Psi \cdot \vc{n}|_{\partial \Omega} = \bfphi \cdot \vc{n}|_{\partial \Omega} = 0
	 \]
	 - the standard potential in the Helmholtz decomposition.
	
	\end{Remark}

\begin{Remark} \label{Rw2}
	
	Unlike $\vr$, $\vr \vu$, and $\vB$, the total entropy $\vr s(\vr, \vt)$ is not weakly continuous in time.
	However, it can be deduced from the entropy inequality that the one sided limits
	\[
	\lim_{t \to \tau - } \intO{ \vr s(\vr, \vt) (t, \cdot) \phi },\
	\lim_{t \to \tau + } \intO{ \vr s(\vr, \vt) (t, \cdot) \phi } 
	\]
	exist for any $\phi \in C(\Ov{\Omega})$, and
	\[
	\lim_{t \to \tau - } \intO{ \vr s(\vr, \vt) (t, \cdot) \phi } \leq \lim_{t \to \tau + } \intO{ \vr s(\vr, \vt) (t, \cdot) \phi }
	\]
whenever $\phi \geq 0$, see \cite[Chapter 3, Section 3.3.3]{FeiNovOpen} for details.
\end{Remark}

\begin{Remark} \label{Rw3}

As we shall see below, the ``magnetic'' ballistic energy $E_{BM}$ introduced in \eqref{w3} is in fact a strictly convex l.s.c function
of the variables $(\vr, S = \vr s(\vr, \vt), \vm = \vr \vu, \vB)$. In particular, its integral is a weakly lower
semi--continuous function; whence \eqref{w14} holds for \emph{any} time $\tau$ with the convention for entropy
discussed Remark \ref{Rw2}.
Alternatively, we may impose a stronger integrated version of \eqref{w14} in the form
\begin{align}
\int_0^T \partial_t \psi	&\intO{ \left( \frac{1}{2} \vr |\vu|^2 + \vr e(\vr, \vt) + \frac{1}{2} |\vB |^2 - \tvt \vr s(\vr, \vt) - \bB \cdot \vB\right) } \dt \br
\quad&
	- \int_0^T \psi \intO{ 	\frac{\tvt}{\vt} \left( \mathbb{S}(\vt, \Grad \vu) : \Grad \vu - \frac{\vc{q}(\vt, \Grad \vt) \cdot \Grad \vt}{\vt} + \zeta (\vt) | \Curl \vB |^2 \right) } \dt \br
	&\quad \geq  \int_0^T \psi \intO{\left( \vr s (\vr, \vt) \partial_t \tvt + \vr s (\vr, \vt) \vu \cdot \Grad \tvt + \frac{\vc{q}(\vr, \vt)}{ \vt} \cdot \Grad \tvt                    \right)      } \dt\br
&\quad + \int_0^\tau \psi \intO{  \Big( \vB \cdot \partial_t \bB - (\vB \times \vu) \cdot \Curl \bB -
	\zeta \Curl \ \vB \cdot \Curl \bB \Big)    } \dt \br 	
	 &\quad -  \int_0^T \psi \intO{ \vr \vc{g} \cdot \vu } \dt \br
	&\quad   - \psi(0) \intO{ \left( \frac{1}{2} \vr_0 |\vu_0|^2 + \vr_0 e(\vr_0, \vt_0) + \frac{1}{2} |\vB_0 |^2 - \tvt (0, \cdot) \vr_0 s
		(\vr_0, \vt_0) - \bB (0, \cdot) \cdot \vB_0  \right)  }
	\label{w15}
\end{align}
for any $\psi \in C^1_c[0,T)$, $\psi \geq 0$.

\end{Remark}

\begin{Remark} \label{Rw4}
	
	Note carefully that the ballistic energy inequality \eqref{w14} remains valid if we replace
	$\bB$ by any other extension $\tvB$ such that
	\[
	\Div \tvB = 0,
	\]
	and
	\[\ \tvB \times \vc{n}|_{\partial \Omega} =
	\bB \times \vc{n}|_{\partial \Omega},\ \mbox{or} \ \tvB \cdot \vc{n}|_{\partial \Omega} =
	\bB \cdot \vc{n}|_{\partial \Omega},
	\]
	respectively. Indeed the difference $\bB - \tvB$ becomes an eligible test function for
	for the weak form of the induction equation \eqref{w11}, endowed with the boundary conditions
	\eqref{w11a}, \eqref{w11b}, respectively.
	
	\end{Remark}

The existence of global--in--time weak solutions in the sense of Definition \ref{De1} will be shown in Section \ref{e} below.

\section{Relative energy}
\label{r}

A suitable form of the \emph{relative energy} for the compressible MHD system reads
\begin{align}
	E &\left( \vr, \vt, \vu, \vc{B} \ \Big| r, \Theta, \vU, \vc{H} \right) =
	\frac{1}{2} \vr |\vu - \vU|^2 + \frac{1}{2} |\vc{B} - \vc{H}|^2 \br &+ 	
	\vr e(\vr, \vt) - \Theta \Big( \vr s(\vr, \vt) - r s(r, \Theta) \Big)
	- \Big( e(r, \Theta) - \Theta s(r, \Theta) + \frac{p(r, \Theta)}{r}     \Big)(\vr - r) - r e(r, \Theta).
	\label{r6}
\end{align}
In applications, the quantity $(\vr, \vt, \vu, \vB)$ stands for a weak solution of the compressible MHD system while
$(r, \Theta, \vU, \vc{H})$ are arbitrary sufficiently smooth functions satisfying the compatibility
conditions
\begin{align}
	r &> 0,\ \Theta > 0 \ \mbox{in}\ [0,T] \times \Ov{\Omega}, \ \Theta|_{\partial \Omega} = \vtB, \ \vU \cdot
	\vc{n}|_{\partial \Omega} = 0, \br
	\Div \vc{H} &= 0 \ \mbox{in}\ (0,T) \times \Omega,\
	\vc{H} \times \vc{n}|_{\partial \Omega} =
	\bB \times \vc{n}|_{\partial \Omega},\  \mbox{or}\ 	\vc{H} \cdot \vc{n}|_{\partial \Omega} =
	\bB \cdot \vc{n}|_{\partial \Omega},\ \mbox{respectively}.
\label{r6a}	
	\end{align}

As a consequence of hypothesis of thermodynamic stability stated in \eqref{c4}, the energy
\[
E(\vr, \vt, \vu, \vB) = \frac{1}{2} \vr |\vu|^2 + \vr e(\vr, \vt) + \frac{1}{2} |\vB|^2
\]
is a strictly convex l.s.c. function if expressed in the variables $(\vr, S = \vr s(\vr, \vt), \vm = \vr \vu, \vB)$, 
see \cite[Chapter 3, Section 3.1]{FeiNovOpen} for details. Moreover, we have
\[
\frac{\partial E(\vr, S, \vm, \vB)}{\partial \vr} = \left(  e(\vr, \vt) - \vt s(\vr, \vt) + \frac{p(\vr, \vt) }{\vr} \right),\ \frac{\partial E(\vr, S, \vm, \vB)}{\partial S} = \vt,
\]
where $\left(  e(\vr, \vt) - \vt s(\vr, \vt) + \frac{p(\vr, \vt) }{\vr} \right)$ is Gibbs' free energy. Consequently, the relative energy can be seen as Bregman divergence associate to the convex function $E$.
In particular,
\[
E \left( \vr, \vt, \vu, \vc{B} \ \Big| r, \Theta, \vU, \vc{H} \right) \geq 0,\ \mbox{and}\
E \left( \vr, \vt, \vu, \vc{B} \ \Big| r, \Theta, \vU, \vc{H} \right)  = 0 \ \Leftrightarrow \
( \vr, \vt, \vu, \vc{B}) = (r, \Theta, \vU, \vc{H} )
\]
as long as $r > 0$.

\subsection{Relative energy inequality}

In the remaining part of this section, we derive a relative energy inequality provided
$(\vr, \vt, \vu, \vB)$ is a weak solution of the compressible MHD system specified in Definition \ref{De1} and $(r, \Theta, \vU, \vc{H})$
arbitrary smooth functions satisfying \eqref{r6a}. Our starting point is
the ballistic energy inequality \eqref{w15}. As noted in Remark \ref{Rw4}, we may flip $\bB$ for $\vc{H}$,
and, of course, $\tvt$ for $\Theta$ obtaining
\begin{align}
	&\intO{ \left( \frac{1}{2} \vr |\vu|^2 + \vr e(\vr, \vt) + \frac{1}{2} |\vB |^2 - \Theta \vr s (\vr, \vt) - \vc{H} \cdot \vB \right) (\tau, \cdot) }
	\br&	+ \int_0^\tau \intO{ 	\frac{\Theta}{\vt} \left( \mathbb{S}(\vt, \Grad \vu) : \Grad \vu - \frac{\vc{q} (\vt, \Grad \vt) \cdot \Grad \vt}{\vt} + \zeta (\vt) | \Curl \vB |^2 \right) } \dt \br
	&\quad \leq - \int_0^\tau \intO{\left( \vr s (\vr, \vt) \partial_t \Theta + \vr s (\vr, \vt) \vu \cdot \Grad \Theta + \frac{\vc{q}(\vt, \Grad \vt)}{ \vt} \cdot \Grad \Theta                    \right)      } \br
	&\quad \quad	 -  \int_0^\tau \intO{  \Big( \vB \cdot \partial_t \vc{H} - (\vB \times \vu) \cdot \Curl \vc{H} -
		\zeta (\vt) \Curl \ \vB \cdot \Curl \vc{H} \Big)    } \dt \br
	&\quad \quad	 + \int_0^\tau\intO{ \vr \vc{g} \cdot \vu } \dt \br
	&\quad \quad + \intO{ \left( \frac{1}{2} \vr_0 |\vu_0|^2 + \vr_0 e(\vr_0, \vt_0) + \frac{1}{2} |\vB_0 |^2 - \Theta (0, \cdot) \vr_0 s
		(\vr_0, \vt_0)  - \vc{H}(0, \cdot) \cdot \vB_0 \right) }
	\label{w14a}
\end{align}

\subsubsection{Momentum perturbation}

Our first goal is to replace $|\vu|^2$ by $|\vu - \vU|^2$. This can be achieved by considering $\vU$ as
a test function in the momentum balance equation \eqref{w9}. After a straightforward manipulation
explained in detail in \cite[Chapter 3, Section 3.2.1]{FeiNovOpen}, we obtain
\begin{align}
	&\intO{ \left( \frac{1}{2} \vr |\vu - \vU|^2 + \vr e(\vr, \vt) + \frac{1}{2} |\vB |^2 - \Theta \vr s (\vr, \vt) -
		\vc{H} \cdot \vc{B} \right) (\tau, \cdot) }
	\br&	+ \int_0^\tau \intO{ 	\frac{\Theta}{\vt} \left( \mathbb{S} (\vt, \Grad \vu) : \Grad \vu - \frac{\vc{q} (\vt, \Grad \vt) \cdot \Grad \vt}{\vt} + \zeta (\vt) | \Curl \vB |^2 \right) } \dt \br
	&\quad \leq - \int_0^\tau \intO{ \Big( \vr (\vu - \vU) \otimes (\vu - \vU)   + p (\vr, \vt) \mathbb{I} - \mathbb{S}(\vt, \Grad \vu) \Big) : \Grad \vU   } \dt \br
	&\quad - \int_0^\tau \intO{ ( \Curl \vB \times \vB ) \cdot \vU  } \dt
	\br
	&\quad- \int_0^\tau \intO{ \vr \Big(  \partial_t \vU + \vU \cdot \Grad \vU \Big) \cdot (\vu - \vU)     } \dt \br
	&\quad - \int_0^\tau \intO{\left( \vr s (\vr, \vt) \partial_t \Theta + \vr s (\vr, \vt) \vu \cdot \Grad \Theta + \frac{\vc{q}(\vt, \Grad \vt)}{ \vt} \cdot \Grad \Theta                    \right)      } + \int_0^\tau\intO{ \vr \vc{g} \cdot (\vu - \vU) } \dt \br
	&\quad  - \int_0^\tau \intO{  \Big( \vB \cdot \partial_t \vc{H} - (\vB \times \vu) \cdot \Curl \vc{H} -
		\zeta (\vt) \Curl \ \vB \cdot \Curl \ \vc{H} \Big)    } \dt \br
	&\quad  + \intO{ \left( \frac{1}{2} \vr_0 |\vu_0 - \vU(0, \cdot) |^2 + \vr_0 e(\vr_0, \vt_0) + \frac{1}{2} |\vB_0 |^2 - \Theta(0, \cdot) \vr_0 s
		(\vr_0, \vt_0) - \vc{H}(0, \cdot) \cdot \vc{B}_0 \right)  }
	\label{r1}
\end{align}
for any $(\Theta, \vU, \vc{H})$ as in \eqref{r6a}.

\subsubsection{Density perturbation}

Next, exactly as in \cite[Chapter 3, Section 3.4]{FeiNovOpen}, we use the quantity
\[
\Big( e(r, \Theta) - \Theta s(r, \Theta) + \frac{p(r, \Theta)}{r}     \Big)
\]
as a test function in the equation of continuity \eqref{w5}:
\begin{align}
	&\intO{ \left( \frac{1}{2} \vr |\vu - \vU|^2 + \vr e (\vr, \vt) + \frac{1}{2} |\vB |^2 - \Theta \vr s (\vr, \vt)
		- \vc{H} \cdot \vc{B} \right) (\tau, \cdot) } \br
	&+ \intO{ \left( \Theta r s(r, \Theta) - \Big( e(r, \Theta) - \Theta s(r, \Theta) + \frac{p(r, \Theta)}{r}     \Big)(\vr - r) - r e(r, \Theta) \right) (\tau, \cdot)     }
	\br&	+ \int_0^\tau \intO{ 	\frac{\Theta}{\vt} \left( \mathbb{S}(\vt, \vu) : \Grad \vu - \frac{\vc{q}(\vt, \Grad \vt) \cdot \Grad \vt}{\vt} + \zeta (\vt) | \Curl \vB |^2 \right) } \dt \br
	&\quad \leq - \int_0^\tau \intO{ \Big( \vr (\vu - \vU) \otimes (\vu - \vU)   + p(\vr, \vt) \mathbb{I} - \mathbb{S}(\vt, \Grad \vu) \Big) : \Grad \vU   } \dt \br
	&\quad + \int_0^\tau \intO{ \frac{\vr}{r} (\vu - \vU) \cdot \Grad p(r, \Theta) } \dt \br
	&\quad - \int_0^\tau \intO{ ( \Curl \vB \times \vB ) \cdot \vU  } \dt
	\br
	&\quad- \int_0^\tau \intO{ \vr \Big(  \partial_t \vU + \vU \cdot \Grad \vU + \frac{1}{r} \Grad p(r, \Theta)\Big) \cdot (\vu - \vU)     } \dt \br
	&\quad - \int_0^\tau \intO{\left( \vr \Big(s(\vr, \vt)  - s(r, \Theta) \Big) \partial_t \Theta + \vr \Big(s (\vr, \vt) - s(r, \Theta) \Big) \vu \cdot \Grad \Theta + \frac{\vc{q}(\vt, \Grad \vt)}{ \vt} \cdot \Grad \Theta                    \right)      }\br &\quad + \int_0^\tau\intO{ \vr \vc{g} \cdot (\vu - \vU) } \dt \br
	&\quad + \int_0^\tau \intO{ \left( \left(1 - \frac{\vr}{r} \right) \partial_t p(r, \Theta) - \frac{\vr}{r} \vu \cdot \Grad p(r,\Theta)      \right)     } \dt
	\br
	&\quad  - \int_0^\tau \intO{  \Big( \vB \cdot \partial_t \vc{H} - (\vB \times \vu) \cdot \Curl \vc{H} -
		\zeta (\vt) \Curl \ \vB \cdot \Curl \ \vc{H} \Big)    } \dt
	\br
	&\quad  + \intO{ \left( \frac{1}{2} \vr_0 |\vu_0 - \vU(0, \cdot) |^2 + \vr_0 e(\vr_0, \vt_0) + \frac{1}{2} |\vB_0 |^2 - \Theta(0, \cdot) \vr_0 s
		(\vr_0, \vt_0) - \vc{H}(0, \cdot) \cdot \vc{B}_0 \right)  } \br
	&\quad 	+ \intO{ \left( \Theta s(r, \Theta) - \Big( e(r, \Theta) - \Theta s(r, \Theta) + \frac{p(r, \Theta)}{r}     \Big)(\vr_0 - r) - r e(r, \Theta) \right) (0, \cdot)     },
	\label{r2}
\end{align}
whenever $r \in C^1([0,T] \times \Ov{\Omega})$, $r > 0$.

\subsubsection{Magnetic field perturbation and the final form of relative energy inequality}

Since
\begin{equation} \label{r4}
 \intO{ \frac{1}{2} |\vc{H}|^2 (\tau, \cdot) } -
 \intO{ \frac{1}{2} |\vc{H}|^2 (0, \cdot) } = \int_0^\tau \intO{ \vc{H} \cdot \partial_t \vc{H} } \dt .
\end{equation}
we may rewrite \eqref{r2} in the final form
\begin{align}
	& \left[ \intO{ E \left( \vr, \vt, \vu, \vc{B} \ \Big| r, \Theta, \vU, \vc{H} \right) } \right]_{t = 0}^{t = \tau}
	\br&	+ \int_0^\tau \intO{ 	\frac{\Theta}{\vt} \left( \mathbb{S}(\vt, \Grad \vu) : \Grad \vu - \frac{\vc{q}
			(\vt, \Grad \vt) \cdot \Grad \vt}{\vt} + \zeta (\vt) | \Curl \vB |^2 \right) } \dt \br
	&\quad \leq - \int_0^\tau \intO{ \Big( \vr (\vu - \vU) \otimes (\vu - \vU)   + p (\vr, \vt) \mathbb{I} - \mathbb{S}(\vt, \Grad \vu) \Big) : \Grad \vU   } \dt \br
	&\quad \quad + \int_0^\tau \intO{ \frac{\vr}{r} (\vu - \vU) \cdot \Grad p(r, \Theta) } \dt \br
	&\quad \quad - \int_0^\tau \intO{ ( \Curl \vB \times \vB ) \cdot \vU  } \dt
	\br
	&\quad \quad- \int_0^\tau \intO{ \vr \Big(  \partial_t \vU + \vU \cdot \Grad \vU + \frac{1}{r} \Grad p(r, \Theta) -
		\vc{g} \Big) \cdot (\vu - \vU)     } \dt \br
	&\quad \quad - \int_0^\tau \intO{\left( \vr \Big( s(\vr, \vt)  - s(r, \Theta) \Big) \partial_t \Theta + \vr \Big(s (\vr, \vt)  - s(r, \Theta) \Big) \vu \cdot \Grad \Theta + \frac{\vc{q} (\vt, \Grad \vt) }{ \vt} \cdot \Grad \Theta                    \right)      } \br
	 	&\quad \quad + \int_0^\tau \intO{ \left( \left(1 - \frac{\vr}{r} \right) \partial_t p(r, \Theta) - \frac{\vr}{r} \vu \cdot \Grad p(r,\Theta)      \right)     } \dt
	\br
	&\quad \quad - \int_0^\tau \intO{  \Big( \vB \cdot \partial_t \vc{H} - (\vB \times \vu) \cdot \Curl \vc{H} -
		\zeta (\vt) \Curl \ \vB \cdot \Curl \ \vc{H} \Big)    } \dt \br &\quad \quad +
	\int_0^\tau \intO{ \vc{H} \cdot \partial_t \vc{H} } \dt \br
		\label{r7}
\end{align}
for any ``test'' functions $(r, \Theta, \vU, \vc{H})$ specified in \eqref{r6a}.

\begin{Remark} \label{Rr1}
	
	In view of future applications, in particular the weak--strong uniqueness principle discussed in the
	forthcoming section, the class of functions  $(r, \Theta, \vU, \vc{H})$ in \eqref{r6a} can be extended
	to the ``maximal regularity'' framework
	\begin{equation} \label{r8}
	r \in W^{1,q}((0,T) \times \Omega), \ (\Theta, \vU, \vB) \in L^q(0,T; W^{2,q}(\Omega; R^7)) \cap W^{1,q}(0,T;
	L^q(\Omega; R^7)),	
		\end{equation}
by a density argument. The exponent $q$ must be taken large enough for all terms in \eqref{r7} to be well defined.	
	\end{Remark}

\section{Weak--strong uniqueness principle}
\label{ws}

The first important property of the weak solutions introduced in Definition \ref{De1} is the \emph{weak--strong uniqueness principle.}
We suppose that $(r, \Theta, \vU, \vc{H})$ is a smooth solution of the problem, specifically,
\begin{align}
\partial_t r + \Div (r \vU) &= 0, \br 	
\partial_t \vU + \vU \cdot \vU + \frac{1}{r} \Grad p(r,\Theta) &=
\frac{1}{r} \Div \mathbb{S}(\Theta, \Grad \vU) + \vc{g} + \frac{1}{r} \Curl \vc{H} \times \vc{H}, \br
\partial_t \vc{H} + \Curl (\vc{H} \times \vU) + \Curl (\zeta (\Theta) \Curl \vc{H} ) &= 0,\
\Div \vc{H} = 0, \br
r \partial_t s(r, \Theta) + r \vc{U} \cdot \Grad s(r, \Theta) &+ \Div \left( \frac{ \vc{q}(\Theta, \Grad \Theta) } {\Theta } \right)\br = \frac{1}{\Theta} \Big( \mathbb{S}(\Theta, \Grad \vU) : \Grad \vU &- \frac{\vc{q} (\Theta, \Grad \Theta \cdot \Grad \Theta) }{\Theta} + \zeta (\Theta) |\Curl \vc{H} |^2 \Big),
	\label{ws1}
	\end{align}
together with the relevant initial and boundary conditions.

\begin{mdframed}[style=MyFrame]
	
	\begin{Theorem} \label{TWS} {\bf (Weak--strong uniqueness)}
		
		Let $\Omega \subset R^3$ be a bounded Lipschitz domain. Suppose the thermodynamic functions
		$p$, $e$, $s$ satisfy the hypotheses \eqref{c1}--\eqref{c6}, and the transport coefficients
		$\mu$, $\eta$, $\kappa$, $\zeta$ comply with \eqref{c8}--\eqref{c11}, with
		\begin{equation} \label{ws25}
			\frac{1}{2} \leq \alpha \leq 1, \ \beta \geq 3.
		\end{equation}
		Let $(\vr, \vt, \vu, \vB)$ be a weak solution of the compressible MHD system \eqref{p1}--\eqref{p4} in $(0,T) \times \Omega$ specified
		in Definition \ref{De1}. Suppose that the same problem, with the same initial and boundary data,
		admits a classical solution $(r, \Theta, \vU, \vc{H})$ in the class
		\begin{equation} \label{ws26}
			r \in C^1([0,T] \times \Ov{\Omega}),\ \Theta \in C^{1,2}([0,T] \times \Ov{\Omega}),\
			\vU, \ \vc{H} \in C^{1,2} ([0,T] \times \Ov{\Omega}; R^3).	
		\end{equation}
		
		Then
		\[
		\vr = r,\ \vu = \vU,\ \vt = \Theta,\ \vc{B} = \vc{H} \ \mbox{a.e. in}\ [0,T] \times \Ov{\Omega}.
		\]

	\end{Theorem}

\end{mdframed}

\bigskip

\begin{Remark} \label{Rwss1}
	
	As we assume existence of a classical solution in Theorem \ref{TWS}, the initial data are also regular,
	in particular,
	\[
	\inf_\Omega \vr_0 > 0,\ \inf_\Omega \vt_0 > 0.
	\]
	
	\end{Remark}

\begin{Remark} \label{Rwss2}
	
	Remarkably, the result holds in the full range of exponents $\alpha$, $\beta$ specified in \eqref{ws25}. 
	In particular, we remove the gap between the weak--strong uniqueness principle stated in \cite{ChauFei}, 
	\cite{FeiNovOpen} for $\beta = 3$ and the existence result requiring $\beta > 6$. The main novelty is 
	formulated in Lemma \ref{AL1} below.
	
	\end{Remark}

The rest of this section is devoted to the proof of Theorem \ref{TWS}. Note that existence
of local in time strong solutions to the compressible MHD system can be established by the nowadays well
understood technique proposed by Valli \cite{VAL}, Valli and Zajaczkowski \cite{VAZA}, or by a more recent
approach via maximal regularity in the spirit of Kotschote \cite{KOT6}. As noted in Remark \ref{Rr1}, the conclusion of Theorem \ref{TWS} remains valid for the strong solutions in the maximal regularity class used in
\cite{KOT6}.

The idea of the proof of Theorem \ref{TWS} is simple: We consider the strong solution $(r, \Theta, \vU, \vc{H})$ as
test functions in the relative energy inequality \eqref{r7} and use Gronwall's type argument. We proceed in several
steps specified below.

\subsection{Magnetic field}

As $\vU$ solves the momentum equation, we get
\begin{align}
- &\int_0^\tau \intO{ \vr \Big(  \partial_t \vU + \vU \cdot \Grad \vU + \frac{1}{r} \Grad p(r, \Theta) -
	\vc{g} \Big) \cdot (\vu - \vU)     } \dt	 \br
&=  \int_0^\tau \intO{ \frac{\vr}{r} \Big( \Div \mathbb{S}(\Theta, \Grad \vU) + \Curl \vc{H} \times \vc{H}   \Big) \cdot (\vU - \vu)     } \dt \br
&= \int_0^\tau \intO{ \left( \frac{\vr}{r} - 1 \right) \Big( \Div \mathbb{S}(\Theta, \Grad \vU) + \Curl \vc{H} \times \vc{H}   \Big) \cdot (\vU - \vu)     } \dt \br
&+ \int_0^\tau \intO{  \Big( \Div \mathbb{S}(\Theta, \Grad \vU) + \Curl \vc{H} \times \vc{H}   \Big) \cdot (\vU - \vu)     } \dt.
\label{ws2}	
	\end{align}

Similarly, $\vc{H}$ being the exact solution of the magnetic field equation,
\begin{align}
 - &\int_0^\tau \intO{  \Big( \vB \cdot \partial_t \vc{H} - (\vB \times \vu) \cdot \Curl \vc{H} -
	\zeta (\vt) \Curl \ \vB \cdot \Curl \ \vc{H} \Big)    } \dt  +
\int_0^\tau \intO{ \vc{H} \cdot \partial_t \vc{H} } \dt	\br
&= \int_0^\tau \intO{ (\vB - \vc{H} ) \cdot \Big( \Curl (\vc{H} \times \vU) + \Curl (\zeta (\Theta) \Curl \vc{H} ) \Big) } \dt \br
&+ \int_0^\tau \intO{ \Big( (\vB \times \vu) \cdot \Curl \vc{H} +
	\zeta (\vt) \Curl \ \vB \cdot \Curl \ \vc{H} \Big) } \dt
	\label{ws3}
	\end{align}

Regrouping the terms in \eqref{r7}, \eqref{ws2}, \eqref{ws3} that do not contain magnetic dissipation, we get
the expression
\begin{align}
- \intO{ (\Curl \vB \times \vB) \cdot \vU } &+ \intO{ ( \Curl \vc{H} \times \vc{H} ) \cdot (\vU - \vu) } \br
&+ \intO{ (\vB - \vc{H}) \cdot \Curl (\vc{H} \times \vU) } + \intO{ (\vB \times\vu) \cdot \Curl \vc{H} } .
\nonumber
	\end{align}
Now, we claim the identity
 \begin{align}
- &\intO{ (\Curl \vB \times \vB) \cdot \vU } + \intO{ ( \Curl \vc{H} \times \vc{H} ) \cdot (\vU - \vu) } \br
&+ \intO{ (\vB - \vc{H}) \cdot \Curl (\vc{H} \times \vU) } + \intO{ (\vB \times\vu) \cdot \Curl \vc{H} } \br
&= \intO{ \Big( \vU \times (\vB - \vc{H} ) \Big) \cdot \Curl ( \vB - \vc{H} )  } +
\intO{ \Big( (\vU - \vu) \times (\vB - \vc{H} ) \Big) \cdot \Curl \vc{H} } \br& +
 \intO{ \Div \left( (\vB - \vc{H} ) \times (\vc{H} \times \vc{U} )\right) }.
\label{ws5}
\end{align}
Indeed the integrals on the right--hand side can be treated as follows:
\begin{align}
&\intO{ \left[ \Big( \vU \times (\vB - \vc{H} ) \Big) \cdot \Curl ( \vB - \vc{H} )   +
 \Big( (\vU - \vu) \times (\vB - \vc{H} ) \Big) \cdot \Curl \vc{H} \right] } \br
&\quad = \intO{ (\vU \times \vB )\cdot \Curl (\vB - \vc{H})} - \intO{ (\vU \times \vc{H} ) \cdot \Curl
	(\vB - \vc{H}) }\br &\quad - \intO{ \Big( (\vU - \vu) \times \Curl \vc{H} \Big) \cdot \vB  } 	
+ \intO{ \Big(  \Curl \vc{H} \times \vc{H} \Big) \cdot (\vU - \vu) } \br
&\quad = - \intO{ (\Curl \vB \times \vB ) \vc{U} } +\intO{ \Big(  \Curl \vc{H} \times \vc{H} \Big) \cdot (\vU - \vu) } \br
&\quad - \intO{(\vU \times \vc{B} ) \cdot \Curl \vc{H} } - \intO{ (\vU \times \vc{H} ) \cdot \Curl
	(\vB - \vc{H}) }\br &\quad - \intO{ \Big( \vU  \times \Curl \vc{H} \Big) \cdot \vB  }  +
\intO{ (\vu \times \Curl \vc{H}) \cdot \vB } \br
&\quad = - \intO{ (\Curl \vB \times \vB )\cdot \vc{U} } +\intO{ \Big(  \Curl \vc{H} \times \vc{H} \Big) \cdot (\vU - \vu) } \br
&\quad  +
\intO{ (\vB \times \vu ) \cdot  \Curl \vc{H} } + \intO{ (\vc{H} \times \vU ) \cdot \Curl
	(\vB - \vc{H}) } \br
&\quad = - \intO{ (\Curl \vB \times \vB )\cdot \vc{U} } +\intO{ \Big(  \Curl \vc{H} \times \vc{H} \Big) \cdot (\vU - \vu) } \br
&\quad  +
\intO{ (\vB \times \vu ) \cdot  \Curl \vc{H} } + \intO{	(\vB - \vc{H}) \cdot \Curl (\vc{H} \times \vU ) } \br
&\quad + \intO{ \Div \left( (\vB - \vc{H} ) \times (\vc{H} \times \vc{U} )\right) }.
\nonumber
		\end{align}

Summing up the previous discussion, we may rewrite inequality \eqref{r7} in the form
\begin{align}
	& \intO{ E \left( \vr, \vt, \vu, \vc{B} \ \Big| r, \Theta, \vU, \vc{H} \right) (\tau, \cdot) }
	\br&	+ \int_0^\tau \intO{ 	\frac{\Theta}{\vt} \left( \mathbb{S} : \Grad \vu - \frac{\vc{q} \cdot \Grad \vt}{\vt} + \zeta (\vt) | \Curl \vB |^2 \right) } \dt \br
	&\quad \leq - \int_0^\tau \intO{ \Big( \vr (\vu - \vU) \otimes (\vu - \vU)   + p \mathbb{I} - \mathbb{S} \Big) : \Grad \vU   } \dt \br
	&\quad \quad + \int_0^\tau \intO{ \frac{\vr}{r} (\vu - \vU) \cdot \Grad p(r, \Theta) } \dt \br
	&\quad \quad- \int_0^\tau \intO{ \Div \mathbb{S}(\Theta, \Grad \vU) \cdot (\vu - \vU)     } \dt \br
	&\quad \quad - \int_0^\tau \intO{\left( \vr (s - s(r, \Theta)) \partial_t \Theta + \vr (s - s(r, \Theta)) \vu \cdot \Grad \Theta + \frac{\vc{q}}{ \vt} \cdot \Grad \Theta                    \right)      } \br
	&\quad \quad + \int_0^\tau \intO{ \left( \left(1 - \frac{\vr}{r} \right) \partial_t p(r, \Theta) - \frac{\vr}{r} \vu \cdot \Grad p(r,\Theta)      \right)     } \dt
	\br
	&\quad \quad + \int_0^\tau \intO{ (\vB - \vc{H} ) \cdot \Curl (\zeta (\Theta) \Curl \vc{H} ) } \dt
	+ \int_0^\tau \intO{ \zeta (\vt) \Curl \vB \cdot \Curl \vc{H} } \dt  \br
	&\quad \quad + \int_0^\tau \intO{ \left( \frac{\vr}{r} - 1 \right) \Big( \Div \mathbb{S}(\Theta, \Grad \vU) + \Curl \vc{H} \times \vc{H}   \Big) \cdot (\vU - \vu)     } \dt \br
	&\quad \quad + \int_0^\tau \intO{ \Big( \vU \times (\vB - \vc{H} ) \Big) \cdot \Curl ( \vB - \vc{H} )  } \dt \br &\quad \quad  +
	\int_0^\tau \intO{ \Big( (\vU - \vu) \times (\vB - \vc{H} ) \Big) \cdot \Curl \vc{H} } \dt \br
	&\quad \quad + \int_0^\tau \intO{ \Div \left( (\vB - \vc{H} ) \times (\vc{H} \times \vc{U} )\right) } \dt
	\label{ws6}
\end{align}

In addition, we have
\begin{align}
&\intO{ (\vB - \vc{H} ) \cdot \Curl (\zeta (\Theta) \Curl \vc{H} ) } =
\intO{ \zeta(\Theta) \Curl \vc{H} \cdot \Curl (\vB - \vc{H}) } \br &\quad +
\intO{ \Div \left( \zeta(\Theta) (\vB - \vc{H}) \cdot \Curl \vc{H} \right) }.
\end{align}
We claim
\begin{align}
\int_0^\tau &\intO{ \Div \left( (\vB - \vc{H} ) \times (\vc{H} \times \vc{U} )\right) } \dt +
\intO{ \Div \left( \zeta(\Theta) (\vB - \vc{H}) \times \Curl \vc{H} \right) }\br	&=
\int_0^\tau \int_{\partial \Omega} \Big( (\vB - \vc{H} ) \times (\vc{H} \times \vc{U} ) \Big) \cdot \vc{n}
+ \Big( \zeta(\Theta) (\vB - \vc{H}) \cdot \Curl \vc{H} \Big) \cdot \vc{n} \ \D \sigma \dt = 0.
\label{ws6a}
	\end{align}
Indeed either $\vB$, $\vc{H}$ satisfy the boundary condition \eqref{p8} and then $(\vB - \vc{H}) \times \vc{n} = 0$ or $\vc{H}$ satisfies the flux condition \eqref{w4a}. In both cases the surface integral in \eqref{ws6a}
vanishes.

Thus relation \eqref{ws6} takes the form
\begin{align}
	& \intO{ E \left( \vr, \vt, \vu, \vc{B} \ \Big| r, \Theta, \vU, \vc{H} \right) (\tau, \cdot) }
	\br&	+ \int_0^\tau \intO{ 	\frac{\Theta}{\vt} \left( \mathbb{S}(\vt, \Grad \vu) : \Grad \vu - \frac{\vc{q}(\vt, \Grad \vt) \cdot \Grad \vt}{\vt} + \zeta (\vt) | \Curl \vB |^2 \right) } \dt \br
	&\quad \quad \leq
	 \int_0^\tau \intO{ \frac{\vr}{r} (\vu - \vU) \cdot \Grad p(r, \Theta) } \dt
	- \int_0^\tau \intO{ p(\vr, \vt) \Div \vU } \dt  \br
	&\quad \quad + \int_0^\tau \intO{ \mathbb{S}(\Theta, \Grad \vU) : \Grad (\vu - \vU)     } \dt
	+ \int_0^\tau \intO{ \mathbb{S}(\vt, \Grad \vu) : \Grad \vU } \dt
	\br
	&\quad \quad - \int_0^\tau \intO{\left( \vr \Big( s(\vr, \vt) - s(r, \Theta) \Big) \partial_t \Theta + \vr \Big(s(\vr, \vt)  - s(r, \Theta) \Big) \vu \cdot \Grad \Theta + \frac{\vc{q}}{ \vt} \cdot \Grad \Theta                    \right)      } \br
	&\quad \quad + \int_0^\tau \intO{ \left( \left(1 - \frac{\vr}{r} \right) \partial_t p(r, \Theta) - \frac{\vr}{r} \vu \cdot \Grad p(r,\Theta)      \right)     } \dt
	\br
	&\quad \quad + \int_0^\tau \intO{ \zeta(\Theta) \Curl (\vB - \vc{H} ) \cdot  \Curl \vc{H}  } \dt
	+ \int_0^\tau \intO{ \zeta (\vt) \Curl \vB \cdot \Curl \vc{H} } \dt  \br
	&\quad \quad + \int_0^\tau \intO{ \left( \frac{\vr}{r} - 1 \right) \Big( \Div \mathbb{S}(\Theta, \Grad \vU) + \Curl \vc{H} \times \vc{H}   \Big) \cdot (\vU - \vu)     } \dt \br
	&\quad \quad + \intO{ \Big( \vU \times (\vB - \vc{H} ) \Big) \cdot \Curl ( \vB - \vc{H} )  } +
	\intO{ \Big( (\vU - \vu) \times (\vB - \vc{H} ) \Big) \cdot \Curl \vc{H} } \br
	&\quad \quad -  \int_0^\tau \intO{ \Big( \vr (\vu - \vU) \otimes (\vu - \vU)   \Big) : \Grad \vU   } \dt
	\label{ws7}
\end{align}

\subsection{Pressure and entropy}

After a simple manipulation, we can rewrite \eqref{ws7} in the form
\begin{align}
	& \intO{ E \left( \vr, \vt, \vu, \vc{B} \ \Big| r, \Theta, \vU, \vc{H} \right) (\tau, \cdot) }
	\br&	+ \int_0^\tau \intO{ 	\frac{\Theta}{\vt} \left( \mathbb{S}(\vt, \Grad \vu) : \Grad \vu - \frac{\vc{q} (\vt, \Grad \vt)\cdot \Grad \vt}{\vt} + \zeta (\vt) | \Curl \vB |^2 \right) } \dt \br
	&\quad \quad \leq
	 \int_0^\tau \intO{ \Big( p(r, \Theta) - p(\vr, \vt) \Big) \Div \vU } \dt  \br
	&\quad \quad - \int_0^\tau \intO{\left( r \Big(s (\vr, \vt) - s(r, \Theta) \Big) \partial_t \Theta + r \Big(s(\vr, \vt) - s(r, \Theta) \Big) \vU \cdot \Grad \Theta                   \right)      } \br
	&\quad \quad + \int_0^\tau \intO{  \left(1 - \frac{\vr}{r} \right) \Big( \partial_t p(r, \Theta) + \vU \cdot \Grad p(r,\Theta)      \Big)     } \dt
	\br
	&\quad \quad + \int_0^\tau \intO{ \mathbb{S}(\Theta, \Grad \vU) : \Grad (\vu - \vU)     } \dt
	+ \int_0^\tau \intO{ \mathbb{S}(\vt, \Grad \vu) : \Grad \vU } \dt
	\br
	&\quad \quad + \int_0^\tau \intO{ \zeta(\Theta) \Curl (\vB - \vc{H} ) \cdot  \Curl \vc{H}  } \dt  + \int_0^\tau \intO{ \zeta (\vt) \Curl \vB  \cdot  \Curl \vc{H}  } \dt \br
	&\quad \quad - \int_0^\tau \intO{ \frac{\vc{q}(\vt, \Grad \vt)}{\vt} \cdot \Grad \Theta } \dt + \int_0^\tau \mathcal{R}_1,
	\label{ws8}
	\end{align}
with a quadratic remainder
\begin{align}
	\mathcal{R}_1 &= \int_0^\tau \intO{ (r - \vr) \Big( s(\vr, \vt)  - s(r, \Theta) \Big) \partial_t \Theta } \dt \br
	&+ \int_0^\tau \intO{ (r - \vr) \Big(s (\vr, \vt) - s(r, \Theta) \Big) \vU \cdot \Grad \Theta } \dt \br
	& + \int_0^\tau \intO{ \vr \Big( s(\vr, \vt) - s(r, \Theta) \Big)(\vU - \vu) \cdot \Grad \Theta } \dt \br
	& + \int_0^\tau \intO{ \left( \frac{\vr}{r} - 1 \right) \Big( \Div \mathbb{S}(\Theta, \Grad \vU) + \Curl \vc{H} \times \vc{H} - \Grad p(r, \Theta)  \Big) \cdot (\vU - \vu)     } \dt \br
	& + \intO{ \Big( \vU \times (\vB - \vc{H} ) \Big) \cdot \Curl ( \vB - \vc{H} )  } +
	\intO{ \Big( (\vU - \vu) \times (\vB - \vc{H} ) \Big) \cdot \Curl \vc{H} } \br
	& -  \int_0^\tau \intO{ \Big( \vr (\vu - \vU) \otimes (\vu - \vU)   \Big) : \Grad \vU   } \dt \br
	& +	\int_0^\tau \intO{ \left( \frac{\vr}{r} - 1 \right) (\vu - \vU) \cdot \Grad p(r, \Theta) } \dt.
	\label{ws9}
\end{align}

Using the identity
\[
\frac{\partial s (r, \Theta) }{\partial \vr} = - \frac{1}{r^2} \frac{\partial p(r, \Theta) }{\partial \vt}
\]
we compute, exactly as in \cite[Section 3.3]{ChauFei},
\begin{align}
	\Big( 1 - &\frac{\vr}{r} \Big) \left( \partial_t p(r, \Theta) + \vU \cdot
	\Grad p(r, \Theta)  \right) + \Div \vU \Big( p(r, \Theta) - p(\vr, \vt) \Big)\br
	= &\Div \vU \Big( p(r, \Theta)  -\frac{\partial p(r, \Theta)}{\partial \vr} (r -\vr ) -\frac{\partial p(r, \Theta)}{\partial \vt} (\Theta -\vt )- p(\vr, \vt) \Big)
	\br
	- &r (r - \vr ) \frac{\partial s (r, \Theta) }{\partial \vr} \Big( \partial_t \Theta + \vU \cdot \Grad \Theta \Big)
	-  r( \Theta - \vt) \frac{\partial s (r, \Theta) }{\partial \vt} \Big( \partial_t \Theta + \vU \cdot \Grad \Theta \Big)
	\br
	- &(\Theta - \vt) \Div \left( \frac{\vc{q}(\Theta, \Grad \Theta) } {\Theta} \right) +
	\left( 1 - \frac{\vt}{\Theta} \right) \left( \mathbb{S}(\Theta, \Grad \vU) : \Grad \vU -
	\frac{ \vc{q} (\Theta, \Grad \Theta)  \cdot \Grad \Theta }{\Theta} \right).
	\nonumber
\end{align}
Consequently, inequality \eqref{ws8} can be written as
\begin{align}
	& \intO{ E \left( \vr, \vt, \vu, \vc{B} \ \Big| r, \Theta, \vU, \vc{H} \right) (\tau, \cdot) }
	\br&	+ \int_0^\tau \intO{ 	\frac{\Theta}{\vt} \left( \mathbb{S}(\vt, \Grad \vu) : \Grad \vu - \frac{\vc{q}
			(\vt, \Grad \vt) \cdot \Grad \vt}{\vt} + \zeta (\vt) | \Curl \vB |^2 \right) } \dt \br
	&\quad \quad \leq
		- \int_0^\tau \intO{ (\Theta - \vt) \Div \left( \frac{\vc{q}(\Theta, \Grad \Theta) } {\Theta} \right) } \dt
		\br
	&\quad \quad + \int_0^\tau \intO{ 	\left( 1 - \frac{\vt}{\Theta} \right) \left( \mathbb{S}(\Theta, \Grad \vU) : \Grad \vU -
		\frac{ \vc{q} (\Theta, \Grad \Theta)  \cdot \Grad \Theta }{\Theta} + \zeta(\Theta) |\Curl \vc{H} |^2 \right)	} \dt \br
	&\quad \quad + \int_0^\tau \intO{ \mathbb{S}(\Theta, \Grad \vU) : \Grad (\vu - \vU)     } \dt
	+ \int_0^\tau \intO{ \mathbb{S}(\vt, \Grad \vu) : \Grad \vU } \dt
	\br
	&\quad \quad + \int_0^\tau \intO{ \zeta(\Theta) \Curl (\vB - \vc{H} ) \cdot  \Curl \vc{H}  } \dt  + \int_0^\tau \intO{ \zeta (\vt) \Curl \vB  \cdot  \Curl \vc{H}  } \dt \br
	&\quad \quad - \int_0^\tau \intO{ \frac{\vc{q}(\vt, \Grad \vt)}{\vt} \cdot \Grad \Theta } \dt + \int_0^\tau \mathcal{R}_2,
	\label{ws10}
\end{align}
with a remainder
 \begin{align}
 	\mathcal{R}_2 &= \int_0^\tau \intO{ (r - \vr) \Big( s(\vr, \vt) - s(r, \Theta) \Big) \partial_t \Theta } \dt \br
 	&+ \int_0^\tau \intO{ (r - \vr) \Big( s(\vr, \vt) - s(r, \Theta) \Big) \vU \cdot \Grad \Theta } \dt \br
 	& + \int_0^\tau \intO{ \vr \Big( s(\vr, \vt) - s(r, \Theta) \Big)(\vU - \vu) \cdot \Grad \Theta } \dt \br
 		& + \int_0^\tau \intO{ \left( \frac{\vr}{r} - 1 \right) \Big( \Div \mathbb{S}(\Theta, \Grad \vU) + \Curl \vc{H} \times \vc{H} - \Grad p(r, \Theta)  \Big) \cdot (\vU - \vu)     } \dt \br
 	& + \intO{ \Big( \vU \times (\vB - \vc{H} ) \Big) \cdot \Curl ( \vB - \vc{H} )  } +
 	\intO{ \Big( (\vU - \vu) \times (\vB - \vc{H} ) \Big) \cdot \Curl \vc{H} } \br
 	& -  \int_0^\tau \intO{ \Big( \vr (\vu - \vU) \otimes (\vu - \vU)   \Big) : \Grad \vU   } \dt \br
 	& +	\int_0^\tau \intO{ \left( \frac{\vr}{r} - 1 \right) (\vu - \vU) \cdot \Grad p(r, \Theta) } \dt \br
 	& + \int_0^\tau \intO{ \Div \vU \Big( p(r, \Theta)  -\frac{\partial p(r, \Theta)}{\partial \vr} (r -\vr ) -\frac{\partial p(r, \Theta)}{\partial \vt} (\Theta -\vt )- p(\vr, \vt) \Big) } \dt \br
 	&+ \int_0^\tau \intO{ r \left( s(r, \Theta) - \frac{\partial s(r, \Theta )}{\partial \vr}(r - \vr) - \frac{\partial s(r, \Theta)}{\partial \vt}(\Theta - \vt) - s(\vr, \vt) \right) (\partial_t \Theta + \vU \cdot \Grad \Theta )                } \dt.
 	\label{ws11}
 \end{align}

\subsection{Diffusive terms}

First, as $\vt$, $\Theta$ coincide on $\partial \Omega$, we can integrate by parts:
\[
- \int_0^\tau \intO{ (\Theta - \vt) \Div \left( \frac{\vc{q}(\Theta, \Grad \Theta) } {\Theta} \right) } \dt =
\int_0^\tau \intO{ \left( \frac{\vc{q}(\Theta, \Grad \Theta) } {\Theta} \right) \cdot (\Grad \Theta - \Grad \vt) } \dt
\]
Putting all diffusive terms on the left--hand side of the relative energy inequality \eqref{ws10}, we get
\begin{align}
\int_0^\tau &\intO{ \left( \frac{\Theta}{\vt} \mathbb{S} (\vt, \Grad \vu) : \Grad \vu  -
 \mathbb{S} (\vt, \Grad \vu) : \Grad \vU)   -
 \mathbb{S}(\Theta, \Grad \vU ) : \Grad (\vu - \vU) \right)  } \dt \br
&= \int_0^\tau \intO{ \left(  \left( \frac{\Theta}{\vt} - 1 \right) \mathbb{S} (\vt, \Grad \vu) : \Grad \vu  + \Big( \mathbb{S} (\vt, \Grad \vu) -  \mathbb{S} (\Theta, \Grad \vU) \Big) : \Grad ( \vu - \vU) \right)    } \dt.
	\label{ws12}	\end{align}
Similarly
\begin{align}
\int_0^\tau &\intO{ \left( \frac{\Theta}{\vt} \zeta(\vt) |\Curl \vB|^2 - \zeta (\Theta)
	\Curl (\vB - \vc{H}) \cdot \Curl \vc{H} - \zeta(\vt) \Curl \vB \cdot \Curl \vc{H}  \right)} \dt \br &=
\int_0^\tau \intO{ \left( \left( \frac{\Theta}{\vt}  - 1 \right) \zeta(\vt) |\Curl \vB|^2
	+ \Big( \zeta(\vt) \Curl \vB - \zeta(\Theta) \Curl \vc{H} \Big) \cdot \Curl ( \vB - \vc{H} ) \right) } \dt.
	\label{ws13}
	\end{align}
Finally,
\begin{align}
- \int_0^\tau &\intO{ \left( \frac{\Theta}{\vt} \frac{\vc{q}(\vr, \Grad \vt) \cdot \Grad \vt }{\vt} \right) } \dt
+  	\int_0^\tau \intO{ (\Theta - \vt) \Div \left( \frac{\vc{q}(\Theta, \Grad \Theta) } {\Theta} \right) } \dt
\br &+ \int_0^\tau \intO{ 	\left( 1 - \frac{\vt}{\Theta} \right) \left(
	\frac{ \vc{q} (\Theta, \Grad \Theta)  \cdot \Grad \Theta }{\Theta} \right)	} \dt +
\int_0^\tau \intO{ \frac{\vc{q}(\vt, \Grad \vt)}{\vt} \cdot \Grad \Theta } \dt \br
&= \int_0^\tau \intO{ \left( \left( 1 - \frac{\Theta}{\vt} \right)  \frac{\vc{q}(\vt, \Grad \vt) \cdot \Grad \vt}{\vt} + \left( 1 - \frac{\vt}{\Theta} \right)
	\frac{ \vc{q} (\Theta, \Grad \Theta)  \cdot \Grad \Theta }{\Theta}  \right)  } \dt \br
&+ \int_0^\tau \intO{ \left( \frac{\vc{q}(\vt, \Grad \vt) }{\vt} - \frac{\vc{q}(\Theta, \Grad \Theta) }{\Theta} \right) \cdot
	(\Grad \Theta - \Grad \vt) } \dt
	\label{ws14}
	\end{align}
Thus \eqref{ws10} can be rewritten in the final form
\begin{align}
	& \intO{ E \left( \vr, \vt, \vu, \vc{B} \ \Big| r, \Theta, \vU, \vc{H} \right) (\tau, \cdot) }
	\br& + \int_0^\tau \intO{ \left( \frac{\Theta}{\vt} - 1 \right) \mathbb{S}(\vt, \Grad \vu) } \dt
	+ \int_0^\tau \intO{ \left( \frac{\vt}{\Theta} - 1 \right) \mathbb{S}(\Theta, \Grad \vU) } \dt \br
	&+ \int_0^\tau \intO{ \left( \frac{\Theta}{\vt} - 1 \right) \zeta(\vt) |\Curl \vB |^2 } \dt
	+ \int_0^\tau \intO{ \left( \frac{\vt}{\Theta} - 1 \right) \zeta(\Theta) |\Curl \vc{H} |^2  } \dt \br
	&+  \int_0^\tau \intO{ \left( \left( 1 - \frac{\Theta}{\vt} \right)  \frac{\vc{q}(\vr, \Grad \vt) \cdot \Grad \vt}{\vt} + \left( 1 - \frac{\vt}{\Theta} \right)
		\frac{ \vc{q} (\Theta, \Grad \Theta)  \cdot \Grad \Theta }{\Theta}  \right)  } \dt  \br
	&+ \int_0^\tau \intO{ \Big( \mathbb{S} (\vt, \Grad \vu) -  \mathbb{S} (\Theta, \Grad \vU) \Big) : \Grad \Big( \vu - \vU \Big)   } \dt  \br
	&+ \int_0^\tau \intO{ \Big( \zeta(\vt) \Curl \vB - \zeta(\Theta) \Curl \vc{H} \Big) \cdot \Curl ( \vB - \vc{H} )  } \dt \br
	&+ \int_0^\tau \intO{ \left( \frac{\vc{q}(\vt, \Grad \vt) }{\vt} - \frac{\vc{q}(\Theta, \Grad \Theta) }{\Theta} \right) \cdot
		(\Grad \Theta - \Grad \vt) } \dt   \br
	& \leq \int_0^\tau \mathcal{R}_2 \dt,
	\label{ws15}
\end{align}
with $\mathcal{R}_2$ given in \eqref{ws11}.

\subsection{Refined estimates}

Similarly to \cite[Chapter 4, Section 4.2.1]{FeiNovOpen}, we introduce the ``essential'' and ``residual'' component
of a measurable function. We set
\begin{align}
\mathcal{O}_{\rm ess} &= \left\{ (t,x) \in (0,T) \times \Omega \ \Big|
 \ \frac{1}{2} \inf_{(0,T) \times \Omega } r \leq \vr(t,x) \leq
 2 \sup_{ (0,T) \times \Omega } r \right. \br
&\left. \frac{1}{2} \inf_{ (0,T) \times \Omega } \Theta  \leq \vt(t,x) \leq
2 \sup_{ (0,T) \times \Omega } \Theta  \right\} \br 	
\mathcal{O}_{\rm res} &= \Big( (0,T) \times \Omega \Big) \setminus \mathcal{O}_{\rm ess},\br
h_{\rm ess} &= h \mathds{1}_{\mathcal{O}_{\rm ess} },\ 	
h_{\rm res} = h - h_{\rm ess} = h \mathds{1}_{\mathcal{O}_{\rm res} }.
\nonumber
	\end{align}

The remaining part of the proof of weak--strong uniqueness principle follows essentially the arguments of
\cite[Chapter 4, Section 4]{FeiNovOpen}, with the necessary modifications to
accommodate the terms containing the magnetic field as well as the complete slip/Dirichlet boundary conditions.
First, we recall the coercivity property of the relative energy:
\begin{align}
E\left(\vr, \vt, \vu, \vB \Big| r , \Theta, \vU, \vc{H} \right)	
&\geq c \Big( | [ \vr - r ]_{\rm ess} |^2 + | [ \vt - \Theta ]_{\rm ess} |^2 +
| [ \vu - \vU ]_{\rm ess} |^2 + | [ \vB - \vc{H} ]_{\rm ess} |^2 \Big) , \br
E\left(\vr, \vt, \vu, \vB \Big| r , \Theta, \vU, \vc{H} \right)
&\geq \left( 1_{\rm res} + [ \vr |\vu |^2 ]_{\rm res} + [ \vr e(\vr, \vt) ]_{\rm res} +
[ \vr |s(\vr, \vt)| ]_{\rm res} + [ | \vB |^2 ]_{\rm res} \right) ,
	\label{ws16}
	\end{align}
where the constant depend on $\inf r$, $\sup r$, $\inf \Theta$, $\sup \Theta$.

\subsubsection{Temperature gradient}

The following crucial estimate was proved in \cite[Chapter 4, Section 4.2.2]{FeiNovOpen} for 
$\beta = 3$. Here, we extend its validity to general $\beta \geq 3$.

\begin{Lemma} \label{AL1}
Under the hypotheses of Theorem \ref{TWS}, there is a constant $c > 0$, 
	\begin{align}
		&\| \vt - \Theta \|_{W^{1,2}(\Omega; R^3)}^2 \br 	
		&\leq  c \intO{ \left( \left( \frac{\Theta}{\vt} - 1 \right)  \frac{\kappa (\vt) |\Grad \vt|^2 }{\vt} + \left(  \frac{\vt}{\Theta} - 1 \right)
			\frac{ \kappa (\Theta) |\Grad \Theta|^2 }{\Theta}  \right)  }  \br
		&+  \intO{ \left(  \frac{\kappa (\Theta) \Grad \Theta }{\Theta} - \frac{\kappa (\vt) \Grad \vt }{\vt}  \right) \cdot
			(\Grad \Theta - \Grad \vt) } + \xi \intO{ E \left( \vr, \vt, \vu, \vB \Big| r, \Theta, \vU, \vc{H} \right)	}
		\label{ws18}	
	\end{align}
	where $\xi > 0$ depends on $\Theta$, $\inf \Theta$, and $\| \vt \|_{L^4(\Omega)}$.
	
\end{Lemma}

\begin{Remark} \label{RRR1}
	
	As we shall see below, the result holds whenever $\kappa (\vt) \approx 1 + \vt^\beta$, $\beta \geq 2$.
	
	\end{Remark}

\begin{proof}
	
	We consider first the essential part. Rewrite  
	\begin{align}
		&\left( \frac{\Theta}{\vt} - 1 \right)  \frac{\kappa (\vt) |\Grad \vt|^2 }{\vt} + \left(  \frac{\vt}{\Theta} - 1 \right)
		\frac{ \kappa (\Theta) |\Grad \Theta|^2 }{\Theta} + 
		\left(  \frac{\kappa (\Theta) \Grad \Theta }{\Theta} - \frac{\kappa (\vt) \Grad \vt }{\vt}  \right) \cdot (\Grad \Theta - \Grad \vt) \br 
		&\geq \Theta \frac{\kappa(\vt)}{\vt^2} |\Grad \vt - \Grad \Theta|^2 + 
		2 \frac{\kappa(\vt)}{\vt^2} (\Theta - \vt) (\Grad \vt - \Grad \Theta) \cdot \Grad \Theta \br 
		&+ \left( \frac{\kappa(\vt) }{\vt} - \frac{\kappa(\Theta) }{\Theta} \right)(\Grad \vt - \Grad \Theta) \cdot \Grad \Theta + \left( \frac{\kappa(\vt) }{\vt} - \frac{\kappa(\Theta) }{\Theta} \right) 
		\left( \frac{\Theta - \vt}{ \Theta } \right) |\Grad \Theta |^2.
		\label{A2} 
	\end{align}	
	Applying H\" older's inequality we deduce from \eqref{A2}
	\begin{align}
		&\intO{  |[ \Grad \vt - \Grad \Theta ]_{\rm ess}|^2 } \br &\aleq \left[ 
		\left( \frac{\Theta}{\vt} - 1 \right)  \frac{\kappa (\vt) |\Grad \vt|^2 }{\vt} + \left(  \frac{\vt}{\Theta} - 1 \right)
		\frac{ \kappa (\Theta) |\Grad \Theta|^2 }{\Theta} + 
		\left(  \frac{\kappa (\Theta) \Grad \Theta }{\Theta} - \frac{\kappa (\vt) \Grad \vt }{\vt}  \right) \cdot (\Grad \Theta - \Grad \vt)	\right]_{\rm ess} \br 
		&+ \xi \intO{ E \left( \vr, \vt, \vu, \vB \Big| r, \Theta, \vU, \vc{H} \right)_{\rm ess}	}, 
		\label{A3}	
	\end{align} 
	where $\xi$ depends on $\Theta$. As a matter of fact, the only hypotheses to be imposed on $\kappa$ is that 
	$\kappa (\vt) > 0$ whenever $\vt > 0$.
	
	As for the residual part, we first rewrite 
	\begin{align}
		&\left( \frac{\Theta}{\vt} - 1 \right)  \frac{\kappa (\vt) |\Grad \vt|^2 }{\vt} + \left(  \frac{\vt}{\Theta} - 1 \right)
		\frac{ \kappa (\Theta) |\Grad \Theta|^2 }{\Theta} + 
		\left(  \frac{\kappa (\Theta) \Grad \Theta }{\Theta} - \frac{\kappa (\vt) \Grad \vt }{\vt}  \right) \cdot (\Grad \Theta - \Grad \vt) \br 
		&= \frac{\Theta}{\vt} \frac{\kappa (\vt) |\Grad \vt|^2 }{\vt} + \frac{\vt}{\Theta}
		\frac{ \kappa (\Theta) |\Grad \Theta|^2 }{\Theta} - \frac{\kappa (\vt) \Grad \vt }{\vt} \cdot \Grad \Theta - 
		\frac{\kappa (\Theta) \Grad \Theta }{\Theta} \cdot \Grad \vt
		\label{A4}
	\end{align}
	Obviously, 
	\[
	\left[ \frac{\vt}{\Theta}
	\frac{ \kappa (\Theta) |\Grad \Theta|^2 }{\Theta} \right]_{\rm res} \aleq 
	\intO{ E \left( \vr, \vt, \vu, \vB \Big| r, \Theta, \vU, \vc{H} \right)_{\rm res}	}.
	\]
	
	Consequently, as $\inf \Theta > 0$, it remains to handle a quantity
	\[
	\left[ \frac{\kappa (\vt) |\Grad \vt|^2 }{\vt^2} - C \kappa (\vt) |\Grad \Theta| \right]_{\rm res}
	\]
	where $C$ is a positive constant depending on $\Theta, \inf \Theta$. Note that this step requires 
	$\vt^2 \aleq \kappa (\vt)$.
	
	To proceed we consider $\Ov{\vt} > 0$ to be fixed below. We write 
	\[
	\left[ \frac{\kappa (\vt) |\Grad \vt|^2 }{\vt^2} - C \kappa (\vt) |\Grad \Theta| \right]_{\rm res} 
	\mathds{1}_{\vt < \Ov{\vt} } + 
	\left[ \frac{\kappa (\vt) |\Grad \vt|^2 }{\vt^2} - C \kappa (\vt) |\Grad \Theta| \right]_{\rm res} 
	\mathds{1}_{\vt \geq \Ov{\vt} } .
	\]
	Obviously, there is $\xi = \xi(\Ov{\vt}, \Theta)$ such that 
	\begin{align} 
		&\intO{ | [ \Grad \vt ]_{\rm res} |^2\mathds{1}_{\vt < \Ov{\vt} } } \br &\aleq \intO{ \left[ \frac{\kappa (\vt) |\Grad \vt|^2 }{\vt^2} - C \kappa (\vt) |\Grad \Theta| \right]_{\rm res} 
			\mathds{1}_{\vt < \Ov{\vt} } } + \xi \intO{ E \left( \vr, \vt, \vu, \vB \Big| r, \Theta, \vU, \vc{H} \right)_{\rm res}	}.
		\label{A5}
	\end{align}
	
	Thus it remains to handle the integral 
	\[
	\intO{ \left[ \frac{\kappa (\vt) |\Grad \vt|^2 }{\vt^2} - C \kappa (\vt) |\Grad \Theta| \right]_{\rm res} 
		\mathds{1}_{\vt \geq \Ov{\vt} } } = 
	\int_{\vt \geq \Ov{\vt} } \left( \frac{\kappa (\vt) |\Grad \vt|^2 }{\vt^2} - C \kappa (\vt) |\Grad \Theta| \right) \dx. 
	\]
	As $\kappa (\vt) \approx (1 + \vt^\beta)$, our task reduces to 
	\[
	\int_{\vt \geq \Ov{\vt} } \left( \frac{\kappa (\vt) |\Grad \vt|^2 }{\vt^2} - C \kappa (\vt) |\Grad \Theta| \right) \dx \approx \int_{\vt \geq \Ov{\vt} } \left( |\Grad \vt^{\frac{\beta}{2}} |^2  - C \vt^\beta |\Grad \Theta| \right) \dx.
	\]
	First observe 
	\[
	\int_{ \vt \geq \Ov{\vt} } \vt^\beta |\Grad \Theta| \dx 
	\leq \int_{ \vt \geq \Ov{\vt} } 
	\left( \vt^{\frac{\beta}{2}} - {\Ov{\vt}}^{\frac{\beta}{2}} \right)^2 |\Grad \Theta| \dx 
	+ \int_{ \vt \geq \Ov{\vt} } \Ov{\vt}^\beta |\Grad \Theta|  \dx, 
	\]
	where 
	\[
	\int_{ \vt \geq \Ov{\vt} } \Ov{\vt}^\beta |\Grad \Theta|  \dx \leq \zeta \intO{ E \left( \vr, \vt, \vu, \vB \Big| r, \Theta, \vU, \vc{H} \right)_{\rm res}	}.
	\]
	
	Finally, it remains to handle  the integral
	\begin{align}
		&\int_{\vt \geq \Ov{\vt} } \left( |\Grad \vt^{\frac{\beta}{2}} |^2  - C \left( \vt^{\frac{\beta}{2}} - {\Ov{\vt}}^{\frac{\beta}{2}} \right)^2|\Grad \Theta| \right) \dx \br &= 
		\int_{\vt \geq \Ov{\vt} } \left( \left| \Grad \left( \vt^{\frac{\beta}{2}} - {\Ov{\vt}}^{\frac{\beta}{2}} \right) \right|^2  - C \left( \vt^{\frac{\beta}{2}} - {\Ov{\vt}}^{\frac{\beta}{2}} \right)^2|\Grad \Theta| \right) \dx.
		\nonumber 	
	\end{align} 
	On the one hand, using Sobolev--Poincar\' e inequality we get
	\begin{align}
		&\int_{\vt \geq \Ov{\vt} }  \left| \Grad \left( \vt^{\frac{\beta}{2}} - {\Ov{\vt}}^{\frac{\beta}{2}} \right) \right|^2 \dx = \intO{  \left| \Grad \left[ \vt^{\frac{\beta}{2}} - {\Ov{\vt}}^{\frac{\beta}{2}} \right]^+ \right|^2                     } \br& \ageq 
		\left( \intO{ \left| \left[ \vt^{\frac{\beta}{2}} - {\Ov{\vt}}^{\frac{\beta}{2}} \right]^+ \right|^6 } \right)^{\frac{1}{3}} = \left( \int_{\vt \geq \Ov{\vt}} \left( \vt^{\frac{\beta}{2}} - {\Ov{\vt}}^{\frac{\beta}{2}} \right)^6   \dx \right)^{\frac{1}{3}}.
		\label{A6}
	\end{align}	
	On the other hand, by H\" older's inequality,
	\[ 
	\int_{\vt \geq \Ov{\vt} } \left( \vt^{\frac{\beta}{2}} - {\Ov{\vt}}^{\frac{\beta}{2}} \right)^2|\Grad \Theta| \dx 
	\leq \left( \int_{\vt \geq \Ov{\vt} } \left( \vt^{\frac{\beta}{2}} - {\Ov{\vt}}^{\frac{\beta}{2}} \right)^4 \dx 
	\right)^{\frac{1}{2}} \left( \int_{\vt \geq \Ov{\vt} } |\Grad \Theta|^2 \dx \right)^{\frac{1}{2}}.  
	\]
	Next, by Jensen's inequality, 
	\begin{align}
		&\left( \int_{\vt \geq \Ov{\vt} } \left( \vt^{\frac{\beta}{2}} - {\Ov{\vt}}^{\frac{\beta}{2}} \right)^4 \dx 
		\right)^{\frac{3}{2}} = 
		\left( \intO{ \left( \left[ \vt^{\frac{\beta}{2}} - {\Ov{\vt}}^{\frac{\beta}{2}} \right]^+ \right)^4 } 
		\right)^{\frac{3}{2}} \br
		&\aleq \intO{ \left( \left[ \vt^{\frac{\beta}{2}} - {\Ov{\vt}}^{\frac{\beta}{2}} \right]^+ \right)^6 } 
		= \int_{\vt \geq \Ov{\vt} } \left( \vt^{\frac{\beta}{2}} - {\Ov{\vt}}^{\frac{\beta}{2}} \right)^6 \dx; 
		\nonumber
	\end{align}
	whence 
	\begin{equation} \label{A7}
		\left( \int_{\vt \geq \Ov{\vt} } \left( \vt^{\frac{\beta}{2}} - {\Ov{\vt}}^{\frac{\beta}{2}} \right)^4 \dx 
		\right)^{\frac{1}{2}} \aleq \left( \int_{\vt \geq \Ov{\vt} } \left( \vt^{\frac{\beta}{2}} - {\Ov{\vt}}^{\frac{\beta}{2}} \right)^6 \dx \right)^{\frac{1}{3}}.
	\end{equation}
	
	Our final observation is that as $\| \vt \|_{L^4(\Omega)}$ is bounded, for any $\delta > 0$, there is 
	$\Ov{\vt}(\delta)$ such that 
	\[
	\int_{\vt \geq \Ov{\vt} } |\Grad \Theta|^2 \dx \leq \delta. 
	\]
	Combining \eqref{A6}, \eqref{A7} we may infer 
	\begin{align}
		&\intO{  | [ \Grad \vt ]_{\rm res} |^2\mathds{1}_{\vt \geq \Ov{\vt} } }	\br &\aleq 
		\intO{ \left[ \frac{\kappa (\vt) |\Grad \vt|^2 }{\vt^2} - C \kappa (\vt) |\Grad \Theta| \right]_{\rm res} 
		} + \xi \intO{ E \left( \vr, \vt, \vu, \vB \Big| r, \Theta, \vU, \vc{H} \right)_{\rm res}	}
		\label{A8}
	\end{align}
	
	Putting together \eqref{A3}, \eqref{A5}, and \eqref{A8} we obtain \eqref{ws18}.
	
\end{proof}

\subsubsection{Viscous stress and magnetic diffusion}

Exactly as in \cite[Chapter 4, Section 4.2.3]{FeiNovOpen}, we estimate the viscous stress dissipation:
\begin{align}
&\left\| \Grad (\vu - \vU) + \Grad^t (\vu - \vU) - \frac{2}{3} \Div (\vu - \vU) \mathbb{I} \right\|_{L^q(\Omega; R^{3 \times 3})}^2 \br	
&\quad \leq c  \left( \intO{ \left( \frac{\Theta}{\vt} - 1 \right) \mathbb{S}(\vt, \Grad \vu)
+ \left( \frac{\vt}{\Theta} - 1 \right) \mathbb{S}(\Theta, \Grad \vU) }  \right. \br
&\quad \left. + \intO{ \Big( \mathbb{S} (\vt, \Grad \vu) -  \mathbb{S} (\Theta, \Grad \vU) \Big) : \Grad \Big( \vu - \vU \Big)   } + \intO{ E\left(\vr, \vt, \vu, \vB \Big| r , \Theta, \vU, \vc{H} \right)	} \right),
	\label{ws19}
	\end{align}
where
\begin{equation} \label{ws20}
	q = \frac{8}{5 - \alpha}.
	\end{equation}

Next, repeating the same arguments, we get
\begin{align}
& \left\|  \Curl ( \vB - \vc{H} ) \right\|^2_{L^2(\Omega; R^3)} \leq c \left(
	\intO{ \left( \frac{\Theta}{\vt} - 1 \right) \zeta(\vt) |\Curl \vB |^2
	+ \left( \frac{\vt}{\Theta} - 1 \right) \zeta(\Theta) |\Curl \vc{H} |^2  } \right. \br
&\quad + \left. \intO{ \Big( \zeta(\vt) \Curl \vB - \zeta(\Theta) \Curl \vc{H} \Big) \cdot \Curl ( \vB - \vc{H} )  } + \intO{ E\left(\vr, \vt, \vu, \vB \Big| r , \Theta, \vU, \vc{H} \right)	} \right) .
\label{ws21}	
	\end{align}

Finally, we recall the generalized Korn--Poincar\' e inequality \cite[Chapter 11, Theorem 11.23]{FeNo6A}:
\begin{equation} \label{ws22}
\left\| \vc{v} \right\|_{W^{1,q}(\Omega; R^3)} \leq c(\delta,q) \left(
\left\| \Grad \vv + \Grad^t \vv - \frac{2}{3} \Div \vv \mathbb{I} \right\|_{L^q(\Omega; R^{3 \times 3})} +
\intO{ \vr |\vc{v}| } \right)
	\end{equation}
whenever
\begin{equation} \label{ws23}
1 < q < \infty, \ \vr \geq 0,\ \intO{ \vr } \geq \delta,\ \intO{ \vr^\gamma } < \frac{1}{\delta}, \ \gamma > 1.
\end{equation}

The previous estimates, together with \eqref{ws23}, allow us to rewrite the inequality \eqref{ws15} in the
form
\begin{align}
	& \intO{ E \left( \vr, \vt, \vu, \vc{B} \ \Big| r, \Theta, \vU, \vc{H} \right) (\tau, \cdot) }
	\br& + \int_0^\tau \left( \left\| \vu - \vU \right\|^2_{W^{1,q}(\Omega; R^3)}  + \left\| \vt - \Theta \right\|^2_{W^{1,2}(\Omega)} + \left\| \Curl ( \vB - \vc{H}) \right\|^2_{L^2(\Omega;R^3)} \right) \dt
	 \br
	& \leq \int_0^\tau \mathcal{R}_2 \dt + c \int_0^\tau  \intO{ E\left(\vr, \vt, \vu, \vB \Big| r , \Theta, \vU, \vc{H} \right)	} dt,\ q = \frac{8}{5 - \alpha},
	\label{ws24}
\end{align} 	
where the remainder $\mathcal{R}_2$ is specified in \eqref{ws11}.

\subsection{Completing the proof of weak--strong uniqueness}

The rest of the proof of weak--strong uniqueness consists in absorbing all integrals appearing in the remainder
$\mathcal{R}_2$ by the left--hand side of \eqref{ws24} and applying Gronwall's argument. We omit the details as the whole procedure is described in \cite[Chapter 4]{FeiNovOpen}.

We have proved Theorem \ref{TWS}.

\begin{Remark} \label{Rws1}
	
	The Lipschitz regularity of the domain $\Omega$ required in Theorem \ref{TWS} may not be sufficient for the
	\emph{existence} of the strong solution to the MHD system.
	
	\end{Remark}

\section{Existence of weak solutions}
\label{e}

The proof of existence of weak solutions can be carried out by means of the approximate scheme, {\it a priori}
estimates, and compactness arguments specified in Chapter 5 and Chapter 12 of the monograph \cite{FeiNovOpen}.

The magnetic field equations can be incorporated exactly as in \cite{DUFE2} without any additional difficulties.
Similarly to \cite[Chapter 12]{FeiNovOpen}, two additional assumptions must be imposed on the constitutive relations, specifically,
\begin{equation} \label{e1}
\lim_{Z \to \infty} \mathcal{S}(Z) = 0,
	\end{equation}
where $\mathcal{S}$ is the function defining the entropy $s_M$. Moreover, we require
\begin{equation} \label{e2}
\beta > 6, 	
	\end{equation}
where $\beta$ is the exponent in \eqref{c10}. While \eqref{e1} corresponds to the Third Law of Thermodynamics,
hypothesis \eqref{e2} is purely technical to ensure the necessary {\it a priori} bounds.

\subsection{Approximation scheme}

The existence of weak solutions can be shown by solving the following multi--level approximation scheme. The equation of continuity \eqref{p1} is replaced by its \emph{parabolic regularization}:
\begin{equation} \label{e3}
	\partial_t \vr + \Div (\vr \vu) = \ep \Del \vr \ \mbox{in}\ (0,T) \times \Omega,\ \ep > 0,
	\end{equation}
supplemented by the Neumann boundary conditions
\begin{equation} \label{e4}
\Grad \vr \cdot \vc{n}|_{\partial \Omega} = 0,
\end{equation}
and the initial condition
\begin{equation} \label{e5}
	\vr(0, \cdot) = \vr_{0, \delta},
	\end{equation}
where  $\vr_{0, \delta} > 0$ is a smooth regularization of the initial density $\vr_0$. Here and hereafter, $\ep$ and $\delta$ are small parameters to be sent to zero successively in the limit passage.

The momentum equation \eqref{p2} is replaced by a Galerkin approximation
\begin{align}
\int_0^\tau &\intO{ \Big[ \vr \vu \cdot \partial_t \bfphi + \vr \vu \otimes \vu : \Grad \bfphi +
	p_\delta (\vr, \vt) \Div \bfphi - \mathbb{S}_\delta (\vt, \Grad \vu) : \Grad \bfphi \Big] } \dt	\br
&- \ep \int_0^\tau \intO{ \Grad \vr \cdot \Grad \vu \cdot \Grad \bfphi } \dt +
\int_0^\tau \intO{ \Big[ \Curl \vc{B} \times \vc{B} \cdot \bfphi + \vr \vc{g} \cdot \bfphi \Big]  } \dt \br
&= \intO{ \vr \vu (\tau, \cdot) \cdot \bfphi } - \intO{ (\vr \vu)_0 \cdot \bfphi(0, \cdot) }
\label{e6}
	\end{align}	
for any test function $\bfphi \in C^1([0,T]; X_N)$, where $X_N$ is a finite dimensional space
spanned by the eigenvectors of the Lam\' e operator with the complete slip boundary conditions \eqref{p5}.
Here,
\begin{equation} \label{e7}
p_\delta (\vr, \vt) = p(\vr, \vt) + \delta \left( \vr^2 + \vr^\Gamma \right),
\mathbb{S}_\delta = \mathbb{S}(\vt, \vu) + \delta \vt \left( \Grad \vu + \Grad^t \vu - \frac{2}{3}
\Div \vu \mathbb{I} \right),\ \delta > 0.	
	\end{equation}

Similarly, the induction equation is replaced by a Galerkin approximation.
We look for solutions in the form
\[
\vc{B} = \bB + \vc{b},
\]
where
\begin{align}
\int_0^\tau &\intO{ \Big[ (\bB + \vc{b}) \cdot \partial_t \bfphi - \Big( (\bB + \vc{b} ) \times \vu \Big) \cdot \Curl \bfphi - \zeta(\vt) \Curl
(\bB + \vc{b} ) \cdot \Curl \bfphi \Big]  } \dt \br
&= \intO{ (\bB + \vc{b} ) \cdot \bfphi (\tau, \cdot) } - \intO{ \vc{B}_0 \cdot \bfphi }
	\label{e8a}
	\end{align}
for any $\bfphi \in C^1([0,T]; Y_N)$.
The finite--dimensional space $Y_N$ is chosen in accordance with the boundary conditions. The functions in
$Y_N$ are smooth, solenoidal,
\[
\Div \bfphi = 0 \ \mbox{in}\ \Omega,
\]
and $\bfphi \times \vc{n}|_{\partial \Omega} = 0$, $\cup_{N > 0} Y_N$ is dense in $H_{0, \tau}$ in the case
of boundary condition \eqref{p8}, and $\bfphi \cdot \vc{n}|_{\partial \Omega} = 0$, $\cup_{N > 0} Y_N$  is dense in $H_{0, n}$ in the case
of boundary conditions \eqref{p9}, \eqref{w4a}.

Finally, the internal energy balance is regularized, adapted to the viscous approximation of the equation of continuity, and solved exactly,
\begin{align}
	\partial_t &(\vr e_{\delta}(\vr, \vt)) + \Div (\vr  e_{\delta}(\vr, \vt) \vu )
	+ \Div \vc{q}_\delta (\vt, \Grad \vt ) = \mathbb{S}_\delta (\vt, \Grad \vu) +
	\zeta(\vt) |\Curl \vB |^2 \br
	&- p(\vr, \vt) \Div \vu + \ep \delta \left( \vr^{\Gamma - 2} + 2 \right) |\Grad \vr|^2 + \delta \frac{1}{\vt^2} - \ep \vt^5,
	 \label{e10}
	\end{align}
with the regularized initial/boundary conditions
\begin{equation} \label{e11}
	\vt(0, \cdot) = \vt_{0, \delta},\ \vt|_{\partial \Omega} = \vtB.
	\end{equation}
Here, we have denoted
\begin{align}
\vc{q}_\delta &= \vc{q} - \delta \left( \vt^\Gamma + \frac{1}{\vt} \right) \Grad \vt , \br
e_\delta &= e + \delta \vt.
\label{e12}	
	\end{align}

\subsection{Existence of weak solutions}

The existence of weak solutions to the Navier--Stokes--Fourier system (without magnetic field) follows by solving
the approximate scheme \eqref{e3}--\eqref{e7}, \eqref{e10}--\eqref{e12}. The proof
in the absence of magnetic field is given with full details in
\cite[Part II, Chapter 5]{FeiNovOpen}. The crucial observation
made in \cite{DUFE2} is that the basic energy estimate consists in considering the scalar product of
\eqref{e6} with $\vu$ and \eqref{e8a} with $\vc{b}$ - a step compatible with the Galerkin approximation;
whence the proof can be reproduced without any essential difficulties in the present setting.
We may therefore state the following result concerning the existence of weak solutions to the compressible MHD system.

\begin{mdframed}[style=MyFrame]
	
	\begin{Theorem} \label{TE} {\bf (Existence of weak solution)}
		
		Let $\Omega \subset R^3$ be a bounded domain of class at least $C^3$. Suppose the  data
		$\vtB$, $\vc{B}_B$ are twice continuously differentiable,
		\[
		{\inf}_{[0,T] \times \partial \Omega} \vtB > 0,\ \Div \bB  = 0 \ \mbox{in}\ (0,T) \times \Omega.
			\]
		Let the thermodynamic functions
		$p$, $e$, $s$ satisfy the hypotheses \eqref{c1}--\eqref{c6} and the transport coefficients
		$\mu$, $\eta$, $\kappa$, and $\zeta$ comply with \eqref{c8}--\eqref{c11}, where
		\begin{equation} \label{e14}
			\frac{1}{2} \leq \alpha \leq 1, \ \beta > 6, \ \lim_{Z \to \infty} \mathcal{S}(Z) = 0.
		\end{equation}
	
	Then for any initial data
	\begin{align}
	\vr_0, \ \vr_0 &> 0,\ \vt_0, \vt_0 > 0,\ \vu_0,\ \vB_0,\ \Div \vB_0 = 0, \br
	&\intO{ \left[ \frac{1}{2} \vr_0 |\vu_0|^2 + \vr_0 e(\vr_0, \vt_0)  +
		\frac{1}{2}|\vB_0 |^2} \right]  < \infty,
	\nonumber
	\end{align}
and any $\vc{g} \in L^\infty((0,T) \times \Omega; R^3)$, the compressible MHD system admits a weak solution in $(0,T) \times \Omega$ in the sense specified in Definition  \ref{De1}. 		
		\end{Theorem}

	\end{mdframed}

\begin{Remark} \label{Re14}
	
	As a matter of fact, positivity of the initial density and temperature is not really necessary in the weak framework. Expressing the initial data in terms of the initial density $\vr_0$, the momentum $\vm_0$,
	the total entropy $S_0$, and the magnetic field $\vB_0$, we can show global in time existence under the condition
	that the initial energy
	\[
	\intO{ \left[ \frac{1}{2} \frac{|\vm_0|^2}{\vr_0} + \vr_0 e(\vr_0, S_0) + \frac{1}{2}|\vB_0|^2 \right] } < \infty
	\]
	is finite.
	\end{Remark}

\begin{Remark} \label{Re15}
	
	The weak solutions resulting from the approximate scheme specified above satisfy the differential version 
	of the ballistic energy inequality \eqref{w15}.
	
\begin{Remark} \label{R16}
	
	The hypothesis $\lim_{Z \to \infty} \mathcal{S}(Z) = 0$ can be omitted and $\beta \geq 3$ allowed as soon as 
	the boundary temperature $\vtB$ is constant.
	
	\end{Remark}

	\end{Remark}

\section{Concluding remarks}
\label{C}

The theory can be extended in a straightforward manner to other types of boundary conditions. 

\begin{itemize}
	
	\item The complete slip boundary conditions \eqref{p5} can be replaced by the Dirichlet (no slip) boundary 
	conditions 
	\[
		\vu|_{\partial \Omega} = \vuB 
		\]
as long as $\vuB \cdot \vc{n} = 0$. This enables to include problems of Taylor--Coutte type. The extension to 
general in/out flow boundary conditions for the velocity would require a more elaborated treatment in view of the 
presence of the magnetic field.	For the in/out flow boundary conditions for the Navier--Stokes--Fourier system in the absence of magnetic field see 
\cite{FeiNovOpen}. 

\item The Dirichlet boundary conditions for the temperature can be replaced by a flux type condition on a part 
of the boundary:
\[
\partial \Omega = \Gamma_1 \cup \Gamma_2,\ \vt|_{\Gamma_1} = \vtB,\ \vc{q} \cdot \vc{n}|_{\Gamma_2} = \vc{h}.
\]
Note that this type of condition is particularly relevant in the modelling of stellar magnetoconvection, 
see Thompson and Christensen--Dalsgaard \cite{ThoChDa}. 

\item The induction flux (electric field) can be prescribed in \eqref{w4a}:
\[
	\Big[  \vB \times \vu + \zeta  \Curl \vB \Big] \times \vc{n}|_{\partial \Omega} = \vc{E}_B \times \vc{n}.	
\]
	
\item Mixed type boundary conditions imposed on \emph{different} components of $\partial \Omega$ can be accommodated without 
essential difficulties.	
	
	\end{itemize}




\def\cprime{$'$} \def\ocirc#1{\ifmmode\setbox0=\hbox{$#1$}\dimen0=\ht0
	\advance\dimen0 by1pt\rlap{\hbox to\wd0{\hss\raise\dimen0
			\hbox{\hskip.2em$\scriptscriptstyle\circ$}\hss}}#1\else {\accent"17 #1}\fi}

\end{document}